\date{}
\title {Conformal restriction: the chordal case}
\author {Gregory Lawler\footnote {Cornell University}
 \and  Oded Schramm\footnote {Microsoft Research}
 \and Wendelin Werner\footnote {Universit\'e Paris-Sud and IUF}}
\documentclass[12pt]{article}
\usepackage{amssymb}
\usepackage{amsmath}
\usepackage{amsthm}
\usepackage{amsfonts}
\usepackage{graphicx}
\input labelfig.tex
\newtheorem{corollary}{Corollary}
\newtheorem{theorem}[corollary]{Theorem}
\newtheorem{lemma}[corollary]{Lemma}
\newtheorem{proposition}[corollary]{Proposition}

\theoremstyle{remark}\newtheorem{definition}[corollary]{Definition}
\theoremstyle{remark}\newtheorem{remark}[corollary]{Remark}
\numberwithin{equation}{section}
\numberwithin{figure}{section}
\numberwithin{corollary}{section}
\newcommand{\Prob} {{\bf P}}
\newcommand{\R}{\mathbb{R}}
\newcommand{\C}{\mathbb{C}}

\newcommand{\U}{\mathbb{U}}
\newcommand {\Disk} {\U}

\newif\ifdraft
\draftfalse
\def\note#1/{\ifdraft {\bf [#1]}\fi}
\def\H {\mathbb{H}}
\def \eps {\epsilon}

\def \P {{\bf P}}

\def \expect {{\bf E}}
\def \E {{\bf E}}
\def \I {{\cal I}}

\def \hulls{{\cal Q}}
\def \p {{\partial}}
\def \diam {{\mathop{\rm diam}}}
\def \dist {{\mathop{\rm dist}}}
\def \Im {{\mathop{\mathrm {Im}}}}
\def \Re {{\mathop{\mathrm {Re}}}}
\def \medb {{\medbreak \noindent}}
\def\GG{G}
\def\aa{\alpha}
\def\ttau{T}
\def \Cp{{\rm cap}}
\def\fillplane{\mathcal F_{\C}}
\def\fillH{\mathcal F_{\H}}
\def\fillHR{\mathcal F^{\R}_{\H}}
\def\fillneg{\mathcal F^{\R_+}_{\H}}
\def \proof {{ \medbreak \noindent {\bf Proof.} }}
\def\proofof#1{{ \medbreak \noindent {\bf Proof of #1.} }}

\def\QED{\qed\medbreak}
\def\eref#1{(\ref {#1})}
\def\SLEkk#1/{$\mathrm{SLE}_{#1}$}
\def\SLEk/{\SLEkk{\kappa}/}
\def\SLE/{$\mathrm{SLE}$}
\def\SLEd#1#2/{$\mathrm{SLE}({#1},{#2})$}
\def\Yl{Y^\lambda} 

\begin{document}
\maketitle

\begin {abstract}
We characterize and describe all 
random subsets $K$ 
of a given simply connected planar domain (the upper half-plane $\H$, say)
which  
satisfy the ``conformal restriction'' property, i.e.,
 $K$ connects two
fixed 
boundary points ($0$ and $\infty$, say) and the law of 
$K$ conditioned to remain in a simply connected
open subset $H$ of $\H$ 
is identical to that of $\Phi(K)$, where $\Phi$ is 
a conformal map from $\H$ onto $H$ with $\Phi(0)=0$ and 
$\Phi(\infty)=\infty$.
The construction of this family 
relies 
on the stochastic Loewner evolution processes with
parameter $\kappa \le 8/3$ and on their distortion 
under conformal maps. We show in 
particular that  
SLE$_{8/3}$ is the only random simple curve satisfying  
conformal restriction and relate it to the outer
boundaries of planar Brownian motion and SLE$_6$.
\end {abstract}



\newpage \tableofcontents \newpage
 
\section {Introduction}

Conformal field theory has been extremely successful in predicting 
the exact values of critical exponents describing the behavior of 
two-dimensional systems from statistical physics.
In particular, in the fundamental papers \cite {BPZ0,BPZ}, 
which were used and extended to the case of the ``surface geometry''
in \cite {Ca1}, it is argued that there is a close relationship between
critical planar systems and some families of conformally 
invariant fields.
This gave rise to 
intense activity both in the theoretical physics community 
(predictions on the exact value of various exponents or 
quantities) and in the
mathematical community (the study of highest-weight representations
of certain Lie algebras).
However, on the mathematical level, the explicit relation 
between the two-dimensional systems and these fields remained
rather mysterious.

More recently, a one-parameter family of random processes called
{\em stochastic Loewner evolution}, or {\em SLE}, was introduced~\cite{S1}. 
The SLE$_\kappa$ process is obtained by solving Loewner's
differential equation with driving term $B(\kappa\, t)$, where
$B$ is one-dimensional Brownian motion, $\kappa>0$.
The SLE processes are continuous, conformally
invariant scaling limits of various discrete curves arising   
in the context of two-dimensional systems.
In particular, for the models studied by physicists for which
conformal field theory (CFT) has been applied and for which 
exponents have been predicted, it is believed that SLE
arises in some way in the scaling limit.  
This has 
been proved
 for  
site-percolation on the triangular lattice~\cite{Sm},
 loop-erased random walks~\cite {LSWlesl}
and  the uniform spanning
tree Peano path~\cite {LSWlesl} (a.k.a.\ the Hamiltonian path on the Manhattan
lattice).  Other models for which this is believed include
the Ising model, the random cluster (or Potts)
models with $q\le 4$, and
the self-avoiding walk.

In a series of papers \cite {LSW1,LSW2,LSW3,LSWa}, the authors 
 derived various properties of the stochastic Loewner evolution 
SLE$_6$, and  used them to compute the 
``intersection exponents'' for   planar Brownian paths.
This program was based on the earlier realization  \cite{LW2}
that any conformally invariant process satisfying a certain
restriction property has  crossing or intersection exponents that are
intimately related to these Brownian intersection exponents.
In particular, \cite {LW2} predicted a strong relation
between planar Brownian motion, self-avoiding walks, and critical
percolation. 
As the boundary of SLE$_6$ is conformally invariant, satisfies
restriction, and can be well understood, computations
of its exponents yielded the Brownian 
intersection exponents
  (in particular, exponents that had been 
predicted by Duplantier-Kwon \cite {DK,Dqg}, 
  disconnection exponents, and
 Mandelbrot's conjecture \cite {Ma}
that the Hausdorff dimension of 
the boundary of planar Brownian motion is $4/3$).
Similarly, the determination of the critical exponents for SLE$_6$
in \cite {LSW1,LSW2,LSW3} combined with Smirnov's~\cite{Sm} proof of 
conformal invariance for critical percolation on the
triangular lattice (along with Kesten's 
hyperscaling relations) facilitated  proofs of several fundamental
properties of critical percolation
\cite {SW, LSW5, S2}, some of which had been predicted in the theoretical 
physics literature
, e.g.,~\cite {dN, N, N2, Pe, SD}.

The main goal of the present paper
 is to investigate  more deeply
 the restriction property that was 
instrumental in relating SLE$_6$ to Brownian motion.
One of our initial motivations was also to  
understand the
scaling limit and  exponents of the two-dimensional self-avoiding walk.
Another motivation was to 
reach a clean understanding of the
relation between SLE and conformal field theory. Consequences
of the present paper in this direction are the subject of
\cite {FW,FW2}. See also \cite {BB,BB2} for aspects of SLE from a CFT
 perspective.  

Let us now briefly describe the conformal restriction property
which we study in the present paper: 
Consider a simply connected domain in the complex plane $\C$, say the 
upper half-plane $\H := \{ x + i y \ : \ y > 0 \}$.
Suppose that two boundary points are given, say $0$ and $\infty$.
We are going to study closed 
random subsets $K$ of $\H$ such that: 
\begin {itemize}
\item
$\overline K \cap \R = \{ 0 \}$, $K$ is unbounded and $\H\setminus K$ has two
connected components.
\item
For all simply connected subsets $H$ of $\H$ such that 
$\H \setminus H$ is bounded and bounded away from the origin,
the law of $K$ conditioned on $K\subset H$ is equal
to the law of $\Phi (K)$, where $\Phi$ is a conformal map 
from $\H$ onto $H$ that preserves the boundary points $0$ and
$\infty$.
\end {itemize}
The law of such a set $K$ is called a (chordal) restriction 
measure.
It turns out that there exists only a one-parameter
 family $\P_\alpha$ of such probability measures, where $\alpha
$ is a positive number, and that 
\begin {equation}
\label {e.chr}
\P_\alpha [ K \subset H 
 ]     
= \Phi'(0)^\alpha
\end {equation}
when 
$\Phi: H \to \H$ is chosen in such a way that
$\Phi(z)/z\to 1 $ as $z\to\infty$.
The measure
$\P_1$  can be constructed easily
by filling in the 
closed loops of a 
Brownian excursion in $\H$, i.e.,
a Brownian motion started from $0$ and conditioned to stay 
in $\H$ for all positive times.

Some of the main results of this paper can be summarized as 
follows:
\begin {enumerate}
\item
The restriction measure $\P_\alpha$ exists if and 
only if $\alpha \ge 5/8$.
\item
The only measure $\P_\alpha$ that is supported on simple
curves is $\P_{5/8}$. It is the law of chordal SLE$_{8/3}$.
\item
The measures $\P_{\alpha}$ for $\alpha > 5/8$ can be 
constructed by adding to the chordal SLE$_\kappa$ curve 
certain Brownian bubbles with intensity $\lambda $, where $\alpha$,
$\lambda$ and $\kappa$ are related by 
$$
\alpha ( \kappa) = \frac {6 - \kappa} {2 \kappa} ,\qquad\qquad
\lambda (\kappa) = \frac {(8 - 3 \kappa)(6 - \kappa)}{2 \kappa}.
$$
\item
For all $\alpha \ge 5/8$, the dimension of the boundary of $K$ 
defined under $\P_\alpha$ is almost surely $4/3$ and locally \lq\lq looks 
like\rq\rq an SLE$_{8/3}$ curve. In particular, the Brownian 
frontier (i.e., the outer boundary of the Brownian path)
looks like a symmetric curve.
\end {enumerate}
As pointed out in \cite {LSWsaw}, this gives strong support to the 
conjecture that chordal SLE$_{8/3}$ is the scaling limit of 
the infinite self-avoiding walk in the upper half-plane and allows one
to recover (modulo this conjecture) the critical exponents
that had been predicted in the theoretical physics 
literature (e.g., \cite {N2, DS}).
This conjecture has recently been tested~\cite{KennedySLEsaw,KennedySAWconf}
by Monte Carlo methods.
Let us also mention (but this will not be the subject of the 
present paper, see \cite {FW,FW2}) that in conformal field theory language,
$-\lambda (\kappa)$ 
is the central charge of the Virasoro algebra  
associated to the 
discrete models (that correspond to SLE$_\kappa$)
and that $\alpha$ is the corresponding highest-weight (for 
a degenerate representation at level 2).

To avoid confusion, let us point out that SLE$_6$ is not a chordal
restriction measure as defined above.  However, it satisfies {\em locality},
which implies a 
 different form of restriction.
We give below a proof of locality for SLE$_6$, which is significantly
simpler than the original proof appearing in~\cite{LSW1}.

We will also study a slightly different restriction property,
which we call {\em right-sided restriction}.  
The measures satisfying right-sided restriction
similarly form a 1-parameter collection $\P_\alpha^+$,
$\alpha>0$.
We present several constructions of the measures
$\P_\alpha^+$.  First, when $\alpha\ge 5/8$, these
can be obtained from the measures $\P_\alpha$
(basically, by keeping only the right-side boundary).
When $\alpha\in (0,1)$, the measure
$\P_\alpha^+$ can also be obtained from an appropriately reflected
Brownian excursion.  It follows that one can
reflect a Brownian excursion off a ray in such a way that
its boundary will have precisely the law of
chordal SLE$_{8/3}$.  A third construction of
$\P_\alpha^+$ (valid for all $\alpha>0$) is
given by a process we call SLE$(8/3,\rho)$.
  The process SLE$(\kappa,\rho)$
is a variant of SLE where a drift is added to the driving
function.  In fact, it is just  Loewner's evolution driven
by a Bessel-type process.
 
The word {\em chordal} refers to connected sets joining two
boundary points of a domain.  There is an analogous
{\em radial} theory, which investigates sets joining an interior
point to the boundary of the domain.  This will be the subject of 
a forthcoming paper \cite {LSWrad}.

We now briefly describe how this paper is organized.
In the preliminary section, we give some definitions, notations
and derive some simple facts that will be used throughout the paper.
In Section~\ref {rest}, we study the family of chordal restriction
measures, and show (\ref {e.chr}).
Section~\ref {bexc} is devoted to the Brownian excursions. We define
these measures and use a result of
B.~Vir\'ag (see~\cite{Virag}) to
show that the filling of such a Brownian excursion
has the law $\P_1$. 

The key to several of the results of the present paper is
the study of the distortion of SLE under conformal
maps, for instance, the evolution of 
the image of the SLE path under the mapping $\Phi$ 
(as long as the SLE path remains in $H$), which is the 
subject of Section~\ref {framesec}.
This study can be considered as a cleaner and more advanced
treatment of similar questions addressed in~\cite{LSW1}.
In particular, we obtain a new short proof of 
the locality property for SLE$_6$, which was essential in
the papers~\cite{LSW1,LSW2,LSW3}.

The SLE distortion behaviour is 
then also used in Section~\ref {s.83} to prove that the law of 
chordal SLE$_{8/3}$ 
is $\P_{5/8}$ and is 
 also instrumental in Section~\ref {s.bub}, where we 
show that all measures $\P_\alpha$ for $\alpha > 5/8$ can be constructed by
adding a Poisson cloud of bubbles to  SLE curves.

The longer Section \ref {s.os} is devoted to the one-sided restriction
measures $\P_\alpha^+$.  
As described above, we
 exhibit various constructions of these measures and show as a by-product
of this description that the two-sided measures $\P_\alpha$ 
do not exist for $\alpha < 5/8$.

A recurring theme in the paper is the principle
that the law $\P$ of a random set $K$ can often be characterized
and understood through the function
$A\mapsto \P[K\cap A\ne\emptyset]$ on an appropriate
collection of sets $A$.  
In Section~\ref {s.fileq} we use this
to show that the outer boundary  of chordal SLE$_6$ is the same
as the outer boundary (frontier) of appropriately
reflected Brownian motion and the outer boundary of
full-plane SLE$_6$ stopped on hitting the unit circle
is the same as the outer boundary of Brownian motion
stopped on hitting the unit circle. 

We conclude the paper with some remarks and pointers to papers 
in preparation.

\section{Preliminaries} 
\label {s.prelim}

In this section
some definitions and notations will be given and 
some basic facts will be recalled. 

\medb
{\bf Important domains.}
The upper half plane $\{x+iy:x\in\R,\,y>0\}$
is denoted by $\H$, the complex plane by
$\C$, the extended complex plane
by $\hat\C=\C\cup\{\infty\}$ and the unit disk by $\U$.

\medb
{\bf Bounded hulls.}
Let $\hulls$ be the set of all bounded 
$A \subset  \overline \H$
such that  $A = \overline {A \cap \H}$
and $\H \setminus
A$ is (connected and) simply connected. 
We call such an $A$ a {\em bounded hull}.

\medb
{\bf The normalized conformal maps $g_A$.}
 For each 
 $A \in \hulls$, there is a 
  unique conformal transformation $g_A:
 \H \setminus A \rightarrow\H$ with $g_A(z)  
- z  \to 0 $ as $z \rightarrow \infty$.  We can 
then define (as in~\cite{LSW1})
\begin{equation}\label{adef}
a(A):= \lim_{z \to \infty} z(g_A(z) - z).
\end{equation}
First note that $a(A)$ is real, because $g_A(z)-z$ has
a power series expansion in $1/z$ near $\infty$ and is real
on the real line in a neighborhood of $\infty$.
Also note that 
\begin{equation}\label{adef2}
a(A)=\lim_{y\to\infty} y\,H(i\,y),
\end{equation}
where $H(z)=\Im \bigl(z-g(z)\bigr)$ is the bounded harmonic function
on $\H\setminus A$ with boundary values $\Im\,z$.
Hence, $a(A)\ge 0$,
and $a(A)$ can be thought of as a measure of the
size of $A$ as seen from infinity.
We will call $a(A)$ the {\em half-plane
capacity} of $A$ (from infinity).
The useful scaling rule for $a(A)$,
\begin{equation}\label{ascal}
a(\lambda\,A)=\lambda^2 \,a(A)
\end{equation}
is easily verified directly.
Since $\Im\, g_A(z)-\Im\, z$ is harmonic, bounded, and has non-positive
boundary values, $\Im\, g_A(z)\le \Im\, z$.  Consequently,
\begin {equation}
\label {e.positive}
0< g_A' (x)  \le 1,\qquad x\in\R\setminus A.
\end {equation}
(In fact, $g_A'(x)$ can be viewed as the probability of
an event, see Proposition \ref{p.P1}.)
 
\medb
{\bf $*$-hulls.}
 Let $\hulls^*$ be the set  of $A \in \hulls$
with $0 \not\in  A$.
We call such an $A$ a $*$-hull.
If $A \in \hulls^*$, then $H = \H \setminus
A$ is as the $H$ in the introduction.

\medb
{\bf The normalized conformal maps $\Phi_A$.}
 For $A \in \hulls^*$, we
define  
$\Phi_A(z) = g_A(z) - g_A(0)$, which
is the unique conformal transformation $\Phi$ of $\H \setminus
A$ onto $\H$ fixing $0$ and $\infty$ with
$\Phi(z) /z\to 1$ as $z \rightarrow \infty$. 

\medb
{\bf Semigroups.}
Let  ${\cal A}$ be the set of all conformal transformations
$\Phi: \H \setminus A \rightarrow \H$ with $\Phi(0) = 0$ 
and $\Phi(\infty) = \infty$, where
$A \in \hulls^*$.  That is,
$\mathcal A=\{\lambda\,\Phi_A:\lambda>0,\,A\in\hulls^*\}$.
Also let
$\mathcal A_1=\{ \Phi_A:A\in\hulls^*\}$.
Note that
${\cal A}$ and ${\cal A}_1$ are both
semigroups under composition.
(Of course, the domain of $\Phi_1\circ \Phi_2$
is $\Phi_2^{-1}(H_1)$ if $H_1$ is the domain of $\Phi_1$.)
  We can consider $\hulls^*$
as a semigroup with the product $\cdot$, where $A \cdot
A'$ is defined by
$\Phi_{A \cdot A'} = \Phi_A \circ \Phi_{A'}$.  
Note that 
\begin{equation}\label{e.aprod}
a(A \cdot A') = a(A) + a(A').
\end{equation}
As $a(A)\ge 0$, this implies that $a(A)$ is monotone in $A$.

\medb
{\bf $\pm$-hulls.}
Let $\hulls_+$ be the set of $A \in \hulls^*$ with
$  A \cap \R  \subset (0,\infty)$.
Let $\sigma$ denote the orthogonal
 reflection about the imaginary
axis, and let 
$\hulls_- =  \{\sigma (A) \ : \  A \in \hulls_+\}$  be
the set of 
$A \in \hulls^*$ with $  A \cap \R \subset (-\infty,0)$.
If $A \in \hulls^*$,  then
we can find unique $A_1, A_3 \in \hulls_+$ and
$A_2,A_4 \in \hulls_-$ such that $A = A_1 \cdot A_2
= A_4 \cdot A_3.$  Note that $\hulls_+,\hulls_-$ are  semigroups.

\medb
{\bf Smooth hulls.}
We will call $A \in \hulls$ a  smooth hull
if there is a smooth curve $\gamma:[0,1] \rightarrow
\C$ with $\gamma(0), \gamma(1) \in \R$,
$\gamma(0,1) \subset \H$,  $\gamma(0,1)$ has no
self-intersections, and 
$\H \cap
\partial A  = \gamma(0,1)$.
Any smooth hull in
$\hulls^*$ is in $\hulls_+ \cup \hulls_-$.
  
\medb
{\bf Fillings.}
If $A\subset\overline\H$ is closed, let $\fillH(A)$
denote the set of all $z\in\overline\H$
such that any path from $z$ to $\infty$
in $\overline\H\cup\{\infty\}$ meets $A$.
In other words, $\fillH(A)$ is the union of $A$
and all the bounded connected components of $\overline\H\setminus A$.
Similarly, $\fillHR(A)$ denotes the union of $A$ with
the connected components of $\overline\H\setminus A$ which
do not intersect $\R$ and $\fillneg(A)$ denotes the union
of $A$ with the connected components of $\overline\H\setminus A$
which do not intersect $[0,\infty)$.
Also, for closed $A\subset \C$,
$\fillplane(A)$ denotes the union of $A$ with the bounded
connected components of $\C\setminus A$.

\medb
{\bf Approximation.}
We
will sometimes 
want to approximate $A \in \hulls$ by smooth
hulls.   The idea of approximating
general domains by smooth hulls is standard
(see, e.g., \cite[Theorem 3.2]{Duren}).

\begin{lemma}  \label{slitlemma1}
Suppose $A \in \hulls_+$.
Then there exists a decreasing
 sequence of smooth hulls
$(A_n)_{n \ge 1}$ such that
$A = \bigcap_{n=1}^\infty A_n $
and 
the increasing sequence $\Phi_{A_n}'(0)$
converges to $\Phi_A'(0)$.
\end{lemma}

\proof The existence of the sequence $A_n$ can be obtained
by various means, for example, by considering the image under
$\Phi_A^{-1}$ of appropriately chosen paths.  
The monotonicity of $\Phi_{A_n}'(0)$ 
follows immediately from the monotonicity of $A_n$
and (\ref {e.positive}).
The convergence is immediate by elementary properties
of conformal maps, since $\Phi_{A_n}$ converges locally uniformly
to $\Phi_A$ on $\H\setminus A$.
\QED

\medb
{\bf Covariant measures.}
Our aim in the present paper is to 
study measures on 
subsets of $\H$. In order to simplify 
further definitions, we give a general definition
that can be applied in various settings. 

Suppose that $\mu$ is a measure on a 
measurable space $\Omega$ whose elements are
subsets of a domain $D$.
Suppose that $\Gamma$ is a set of conformal transformations
from subdomains  $D' \subset D$ onto $D$ that is closed 
under composition. 
We say that $\mu$ is {\em covariant under $\Gamma$}
(or {\em $\Gamma$-covariant}) if for all $\varphi \in
\Gamma$, 
the measure $\mu$ restricted to the 
set  $\varphi^{-1} (\Omega) := \{\varphi^{-1}(K) : K \subset D\}$
is equal  to a constant
 $F_{\varphi}$ times the image measure 
$\mu \circ \varphi^{-1}$. 

If $\mu$ is a finite $\Gamma$-covariant measure, then
$F_\varphi = \mu[\varphi^{-1}(\Omega)]/
\mu[\Omega]$.  Note that a probability measure $\P$ on $\Omega$
 is $\Gamma$-covariant
 if and only if for all $\varphi \in \Gamma$
with $F_\varphi = \P [ \varphi^{-1} (\Omega) ] >0$, 
the conditional law of $\P$ on $\varphi^{-1} (\Omega)$ is 
equal to $\P \circ \varphi^{-1}$.  

Also note that if $\mu$ is covariant under $
\Gamma$, then
$ F_{\varphi \circ \psi} = F_{\varphi} F_{\psi}$
for all $\varphi, \psi \in \Gamma$, 
because the 
image measure of $\mu$ under 
$\varphi^{-1}$
is  
$ F_\varphi^{-1} \, \mu$ restricted to $\varphi^{-1} (\Omega)$, 
so that the image under $\psi^{-1}$ of this measure
is $F_\varphi^{-1}\,  F_\psi^{-1}\, \mu$
restricted to $\psi^{-1} \circ \varphi^{-1} (\Omega)$.
Hence, the mapping $F : \varphi \mapsto F_\varphi$ is 
a semigroup homomorphism from $\Gamma$
 into the commutative 
multiplicative semigroup $[0, \infty)$.
When $\mu$ is a probability measure,  this mapping is into $[0,1]$.

We say that a measure $\mu$ is {\em $\Gamma$-invariant}
if it is $\Gamma$-covariant with $F_\varphi \equiv 1$.

\medb
{\bf Chordal Loewner chains.}
Throughout this paper, we will make use of chordal Loewner chains.
Let us very briefly recall their definition (see \cite {LSW1} for details).
Suppose that $W=(W_t, t \ge 0)$ is a real-valued continuous function.
Define for each $z \in \overline \H$, the solution $g_t (z)$ of the initial 
value problem
\begin {equation}
\label {e.gt}
\p_t g_t (z) = \frac {2} {g_t (z) - W_t}, \qquad g_0(z) =z 
.\end {equation}
For each $z\in\overline\H$ there
is a time $\tau=\tau(z)\in[0,\infty]$ such that the
solution $g_t(z)$ exists for $t\in[0,\tau]$ and
$\lim_{t\nearrow\tau}g_t(z)=W_\tau$ if $\tau<\infty$. 
The evolving hull of the Loewner evolution is defined 
as $K_t:=\{z\in\overline\H:\tau(z)\le t\}$, $t\ge 0$. 
It is not hard to check that $K_t\in\hulls$.

Then, it is easy to see that $g_t$ is the unique conformal map from 
$\H \setminus K_t$ onto $\H$ such that $g_t (z) = z + o(1)$  
when $z \to \infty$.  
Also, $a(K_t) = 2t$. 
When the function $W$ is chosen to be $W_t = \sqrt {\kappa} B_t$, 
 where $B$ is 
a standard one-dimensional Brownian motion, then the corresponding random Loewner
chain is chordal SLE$_\kappa$ (SLE stands for stochastic Loewner evolution).

\section {Two-sided restriction}
\label {rest}

In this section we will be studying certain probability measures
on a collection $\Omega$ of subsets of $\H$.  We start by defining
$\Omega$.

\begin{definition}\label{omegadef}
Let $\Omega$ be the collection of relatively closed subsets
$K$ of\/ $\H$ such that 
\begin {enumerate}
\item 
$K$ is connected, 
$\overline K \cap \R = \{ 0 \}$ and $K$ is unbounded.
\item
$\C\setminus \overline K$ is connected.
\end {enumerate}
\end{definition}

A simple example of a set $K \in \Omega$ is a simple 
curve $\gamma$ from $0$ to infinity in the upper half-plane.
If $\gamma$ is just a curve from zero to infinity in the 
upper half-plane with double-points, then one can 
take $K=\fillHR(\gamma)\in\Omega$, which is
the set obtained by filling in the loops
created by $\gamma$.

We endow $\Omega$ with the $\sigma$-field generated by the 
events 
$\{K \in \Omega \ : \ 
 K \cap A = \emptyset \}$, where $A \in \hulls^*$.
It is easy to check that this family of events is closed under
finite intersection, so that a probability measure on $\Omega$ 
is characterized by the values of 
$\P [ K \cap A = \emptyset ]$ for $A \in \hulls^*$.
Thus: 

\begin{lemma}\label{l.uniq}
Let $\P$ and $\P'$ be two probability measures
on $\Omega$.  If $\P[K\cap A=\emptyset]=\P'[K\cap A=\emptyset]$
holds for every $A\in\hulls^*$, then $\P=\P'$. \QED
\end{lemma}

It is worthwhile to note that the $\sigma$-field on $\Omega$ is
the same as the Borel $\sigma$-field induced by the Hausdorff metric
on closed subsets of $\overline\H\cup\{\infty\}$.

\begin {proposition}
\label {p.first} 
Let $\Gamma$ be the semigroup of dilations, $z \mapsto \lambda z,
\lambda > 0$.
For any probability measure $\P$ on $\Omega$,
the following four statements are equivalent:
\begin {enumerate}
\item $\P$ is  $\Gamma$-invariant and ${\cal A}_1$-covariant.
\item $\P$ is ${\cal A}$-covariant. 
\item  
There exists an $\alpha  > 0$
such that for  
all $A \in \hulls^*$,
$$
\P [ K \cap A = \emptyset ] = \Phi_A'(0)^{\alpha}
.
$$
\item
 There exists an $\alpha > 0$ such that for all smooth 
hulls $A\in
  {\hulls}^*$,
$$
\P [ K \cap A = \emptyset ] = \Phi_A'(0)^{\alpha }.
$$
\end {enumerate}
Moreover, for each fixed $\alpha > 0$, there exists at most one 
probability measure $\P_{\alpha }$ satisfying these conditions.
\end {proposition}

\begin{definition}
If the measure $\P_\alpha$ exists, we call it
the {\em two-sided restriction measure with exponent $\alpha$}.
\end{definition}

\proof
Lemma~\ref{l.uniq} shows that a measure satisfying 3 is unique.
A probability measure is $\Gamma$-covariant if
and only if it is $\Gamma$-invariant.  Therefore, 1 and 2 are equivalent.
As noted above, any $A\in\hulls^*$ can be written as $A_+\cdot A_-$
with $A_\pm\in\hulls_\pm$.  Using this and Lemma~\ref{slitlemma1}, we
may deduce that conditions
 3 and 4 are also equivalent.
Since $\Phi_{\lambda A}'(0)= \Phi_A'(0)$ for $A\in\hulls^*$, $\lambda>0$,
3 together with Lemma~\ref{l.uniq} imply
that $\P$ is $\Gamma$-invariant. Because
$\Phi_{A_1\cdot A_2}'(0) = \Phi_{A_1}'(0)\, \Phi_{A_2}'(0)$, 
3 also implies that for all $A_1,A_2\in\hulls^*$,
\[     \Prob[K \cap (A_1\cdot A_2) = \emptyset]
           = \Prob[K\cap A_1 = \emptyset ]\, \Prob[K 
    \cap A_2 = \emptyset], \]
which implies 1.
   Hence, it suffices to show
that 1 implies 4.

Suppose 1 holds.  Define the homomorphism $F$ of ${\hulls}^*$
onto the multiplicative semigroup $(0,1]$ 
by $F( A) = \P [ K \cap  A = \emptyset]$.
We also write $F(\Phi_A)$ for $F(A)$.
Let $G_t (z)$ be the solution of the initial value problem 
\begin {equation}\label{e.Geq} 
\p_t G_t (z) = \frac {2 G_t (z)} {G_t (z) - 1} , \qquad G_0 (z) = z,
\end {equation}
for $z \in \overline \H$.
Note that this function can equivalently be defined as 
$G_t(z)=g_t(z)-g_t(0)=g_t(z)+2 t$, where $(g_t)$ is  
the chordal Loewner chain driven by the function $W_t = 1 -2t$. 
Hence, $G_t$ is the unique conformal 
map from $\H \setminus K_t$ onto $\H$ such that 
$G_t (0)=0$ and $G_t(z) /z\to 1$ when $z \to \infty$. 
(Here, $K_t$ is the evolving hull of $g_t$.) 
Also, and this is why we focus on these functions $G_t$, one has 
$G_t\circ G_s=G_{t+s}$ in $\H \setminus K_{t+s}$, for all $s,t\ge 0$.
Since $F$ is a homomorphism,
this implies that
$F(G_t)= \exp(-2\,\alpha\, t)$ for some constant $\alpha\ge 0$
and all $t\ge 0$, or that $F(G_t)=0$ for all $t>0$.
However, the latter possibility would imply that
$K\cap K_t\ne\emptyset$ a.s., for all $t>0$. 
Since $\bigcap_{t>0}K_t=\{1\}$ and $1\notin\overline K$, 
this is ruled out.  Hence,
$F(G_t)=\exp(-2\,\alpha\, t)$, $t\ge 0$.
Differentiating~\eref{e.Geq} with respect to $z$ gives  
$G_t'(0)=\exp(-2\,t)$.  Thus $F(G_t)=G_t'(0)^\alpha$. 

Now, set $G_t^\lambda(z)=\lambda\,G_t(\lambda^{-1} z)$, $\lambda>0$.
Then $G_t^\lambda:\H\setminus\lambda K_t\to\H$ is a suitably normalized 
conformal map.  By our assumption of scale invariance of the
law of $K$, we have
$F(G_t^\lambda)=F(G_t)=G_t'(0)^\alpha=({G_t^\lambda})'(0)^\alpha$.
We may therefore conclude that
$F(\Phi_A)=\Phi_A'(0)^\alpha$ for every $A$ in
the semigroup $\mathcal A_0$
generated by $\{\lambda K_t:t\ge 0,\, \lambda>0\}$. 
To deduce that  
\begin{equation}\label{bigspace}
\forall A\in \hulls_+\,,\qquad
F(\Phi_A)=\Phi_A'(0)^\alpha,
\end{equation}
we rely on the following lemma:

\begin {lemma}
\label {l.suppl}
There exists 
a topology on $\hulls_+$
for which $\mathcal A_0$ 
is dense, $F$ is continuous, and
$\Phi_A\mapsto \Phi_A'(0)$ is continuous.
\end {lemma}

\noindent
{\bf Proof of Lemma \ref{l.suppl}.}
Given $A\in\hulls_+$ and a sequence $\{A_n\}\subset\hulls_+$,
we say that $A_n\in\hulls_+$ converges to $A\in\hulls_+$
if $\Phi_{A_n}$
converges to $\Phi_A$ uniformly on compact subsets of
$\overline\H\setminus A$
and $\bigcup_n A_n$ is bounded away from $0$ and $\infty$.
(This is very closely related to what is known as
the Carath\'eodory topology.)

Now assume that $A_n\to A$, where $A_n,A\in\hulls_+$.
It is immediate that $\Phi_{A_n}'(0)\to\Phi_A'(0)$, by Cauchy's
derivative formula (the maps may be extended to a neighborhood
of $0$ by Schwarz reflection in the real line).
Set $A^+_n= \Phi_{A_n}(A\setminus A_n)$ and
$A^-_n= \Phi_A(A_n\setminus A)$.
We claim that there is a constant $\delta>0$
and a sequence $\delta_n\to 0$ such that 
\begin{equation}\label{e.AA}
A^+_n\cup A^-_n\subset
\{x+iy:x\in[\delta,1/\delta],\,y\le \delta_n\}
\,.
\end{equation}
Indeed, since the map $\Phi_{A_n}\circ \Phi_A^{-1}$
converges to the identity, locally uniformly in $\H$,
it follows (e.g., from the argument principle) that
for every compact set $S\subset\H$ for all sufficiently large
$n$, $S$ is contained in the image of $\Phi_{A_n}\circ \Phi_A^{-1}$,
which means that $A^+_n\cap S=\emptyset$.
Similarly, $\Phi_{A}^{-1}\circ\Phi_{A_n}$ converges locally
uniformly in $\overline\H\setminus A$ to the identity map,
and this implies that $A^-_n\cap S=\emptyset$ for all
sufficiently large $n$.
It is easy to verify that $A^+_n\cup A^-_n$ is bounded and
bounded away from $0$.  Consequently, we have~\eref{e.AA}
for some fixed $\delta>0$ and some sequence $\delta_n\to 0$.

Suppose $\limsup_{n\to\infty} F(A_n)>F(A)$.
Then $\limsup_{n\to\infty} \P[K\cap A_n\ne\emptyset = K\cap A]>0$.
By mapping over with
$\Phi_A$ and using 1, it then follows that there is some
$\eps>0$ such that for infinitely many $n$
$\P[K\cap A_n^-]>\eps$. 
Therefore, with positive probability, 
$K$ intersects infinitely many $A_n^-$. 
Since $K$ is closed and~\eref{e.AA} holds, this would then imply that
$\P\bigl[K\cap[\delta,1/\delta]\bigr]>0$, a contradiction.
Thus $\limsup_{n\to\infty} F(A_n)\le F(A)$.  A similar argument
also shows that
$\liminf_{n\to\infty} F(A_n)\ge F(A)$, and so
$\lim_{n\to\infty} F(A_n)=F(A)$, and the continuity of $F$ is
verified.

To complete the proof of the Lemma, we now show that
$\mathcal A_0$ is dense in $\hulls_+$.  Let $A\in\hulls_+$.
Set $A':=A\cup [x_0,x_1]$, where $x_0:=\inf (A\cap \R)$
and $x_1:=\sup (A\cap \R)$.  For $\delta>0$, $\delta<\Phi_A(x_0)/2$,
let $D_\delta$ be the set of points in $\H$ with
distance at most $\delta$ from $[\Phi_A(x_0),\Phi_A(x_1)]$.
Let $E_\delta$ denote the closure of $A\cup\Phi_A^{-1}(D_\delta)$.
It is clear that $E_\delta\to A$ as $\delta\to 0+$ in the topology
considered above.  
It thus suffices to approximate $E_\delta$.
Note that $\overline {\p E_\delta\cap\H}$ is a simple path,
say $\beta:[0,s]\to\overline\H$ with $\beta(0),\beta(s)\in\R$.
We may assume that $\beta$ is parametrized by half-plane capacity
from $\infty$,
so that $a\bigl(\beta[0,t]\bigr)=2\,t$, $t\in[0,s]$.
Set $g_t:=g_{\beta[0,t]}$,
$\Phi_t:=\Phi_{\beta[0,t]}=g_t-g_t(0)$, 
$U_t:= g_t\bigl(\beta(t)\bigr)$,
$\tilde U_t:=U_t-g_t(0)=\Phi_t\bigl(\beta(t)\bigr)$, $t\in[0,s]$.
By the chordal version of Loewner's theorem, we have
$$
\p_t g_{t}(z)= \frac 2{g_{t}(z)-U_t}\,.
$$
Thus,
\begin{equation}\label{e.lu}
\p_t \Phi_t(z)=
\frac 2{\Phi_{t}(z)-\tilde U_t}+
\frac 2{\tilde U_t}= 
\frac {2\, \Phi_t(z)}{(\Phi_t(z)-\tilde U_t)\,\tilde U_t},\qquad
\Phi_0(z)=z
\,.
\end{equation}
Since $\tilde U_t$ is continuous and positive, there is a sequence of
piecewise constant functions
$\tilde U^{(n)}:[0,s]\to(0,\infty)$ 
such that $\sup\{|\tilde U^{(n)}_t-\tilde U_t|:t\in[0,s]\}\to 0$
as $n\to\infty$.  Let $\Phi^{(n)}_t$ be the solution of~\eref{e.lu} with
$\tilde U^{(n)}_t$ replacing $\tilde U_t$.  
Then, clearly, $\Phi^{(n)}_s(z)\to \Phi_s(z)=\Phi_{E_\delta}$
locally uniformly in $\overline \H\setminus{E_\delta}$.
Note that the solution of~\eref{e.lu} with $\tilde U_t$ constant is of
the form $G^\lambda_{t'}$, where $\lambda=\tilde U_0$ and $t'$ is some function of
$t$ and $\tilde U_0$.  It follows that $\Phi^{(n)}_s$ is in the semigroup generated
by $\{G^\lambda_t:\lambda>0,t\ge 0\}$.  Hence, $\mathcal A_0$
is dense in $\hulls_+$ and Lemma \ref {l.suppl} is established.
\qed

\medbreak
\noindent
{\bf End of the proof of Proposition \ref {p.first}.}
Clearly, the lemma implies~\eref{bigspace}.
By symmetry, there is a constant $\alpha_-$
such that $F(A)=\Phi_A'(0)^{\alpha_-}$ holds for every
$A\in\hulls_-$.
To verify that $\alpha_-=\alpha$,
 let $\epsilon$ be small, let
$A_+=\{e^{i\theta}: \theta\in [0,\pi/2-\epsilon]\}$,
$A_-=\{e^{i\theta}: \theta\in [\pi/2+\epsilon,\pi]\}$
and $A=A_+\cup A_-$.
Set $A_+^*=\Phi_{A_-}(A_+)\in\hulls_+$
and $A_-^*=\Phi_{A_+}(A_-)\in\hulls_-$.
Note that $\Phi_{A_+}'(0)$ and $\Phi_{A_+}'(0)$ are
bounded away from zero, but $\lim_{\eps\searrow 0}\Phi_A'(0)= 0$.
As 
\begin{equation}\label{phicomp}
\Phi_{A_-^*}\circ \Phi_{A_+}=\Phi_A= \Phi_{A_+^*}\circ \Phi_{A_-}\,,
\end{equation}
we have $\Phi_{A_-^*}'(0)\to 0$ 
when $\eps\to 0$.
By applying $F$ to~\eref{phicomp} we get
$$
\Phi_{A_-^*}'(0)^{\alpha_-}\;\Phi_{A_+}'(0)^\alpha=
\Phi_{A_+^*}'(0)^\alpha
\;\Phi_{A_-}'(0)^{\alpha_-}\,.
$$
As $\Phi_{A_-^*}'(0)= \Phi_{A_+^*}'(0)$ (by symmetry), this
means that $\Phi_{A_-^*}'(0)^{\alpha-\alpha_-}$
stays bounded and bounded away from zero as $\eps\searrow 0$, which gives
 $\alpha_-=\alpha$.
Since every $A\in\hulls^*$ can be written as
$A_+\cdot A_-$ this establishes 4 with $\alpha\ge0$.
The case $\alpha=0$ clearly implies $K=\emptyset$ a.s., which
is not permitted.  This completes the proof.
\QED

Let us now conclude this section with some simple remarks:
\begin {remark}
If $K_1, \ldots , K_n$ 
are independent sets with respective laws
 $\P_{\alpha_1}$,
\ldots,$\P_{\alpha_n}$, then the law of  
the filling  $K:= 
\fillHR (K_1 \cup \ldots \cup  
K_n)$ of the union of the $K_j$'s
 is $\P_{\alpha}$ with $\alpha= \alpha_1+\cdots + \alpha_n$
because
$$
\P [  K \subset \Phi^{-1} (\H) ]
= \prod_{j=1}^n \P [ K_j \subset \Phi^{-1} (\H) ] 
= \Phi'(0)^{\alpha_1 + \cdots + \alpha_n}.
$$
\end {remark}

\begin {remark}\label{r.nolesshalf}
When $\alpha < 1/2$, 
the measure $\P_\alpha$ does not exist.
To see this, suppose it did.  Since 
  it is unique, it is
invariant under the symmetry $\sigma: x+iy \mapsto -x+iy$.
Let $A =  \{ e^{i \theta} \ : \ \theta \in [0,\pi/2] \} $.
Since $K$ is almost surely connected and joins $0$ to infinity,
it meets either  $A$ or $\sigma(A)$. Hence, symmetry implies that
$$ 
\Phi' _A(0)^\alpha = 
\P [ K \cap A = \emptyset] \le 1/2.$$
On the other hand, one can calculate directly
$\Phi'_A(0) = 1/4$, and hence 
$
\alpha \ge 1/2
$.
We will show later in the paper (Corollary
\ref{mar26.cor1}) that $\P_\alpha$ 
only exists for $\alpha \ge 5/8$.
\end {remark}

\begin {remark}
We have chosen to study subsets of the upper half-plane with
the two special boundary points $0$ and 
$\infty$, but our analysis clearly 
applies to any simply connected domain $O \not= \C$
with two distinguished boundary points $a$ and $b$, $a\ne b$.
(We need to assume that the boundary of $O$ is sufficiently nice near
$a$ and $b$.  Otherwise, one needs to discuss prime ends in place of
the distinguished points.)
For instance, if $\partial O$ is smooth in the neighborhood of 
$a$ and $b$, then for a conformal map $\Phi$ from 
a subset $O'$ of $O$ onto $O$, we get 
$$
\P  [ K \cap (O \setminus O') = \emptyset 
] = 
( \Phi'(a) \Phi'(b))^\alpha\,,
$$
where $\P$ denotes the image of $\P_\alpha$ under a conformal
map from $\H$ to $O$ that takes $a$ to $0$ and $b$ to $\infty$.
\end {remark}

\begin {remark}
The proof actually shows that  weaker assumptions
on $\Omega$ are sufficient for the proposition.
Define $\Omega^b$ just as $\Omega$ was defined, except that Condition 1 is 
replaced by the requirements that $K\ne\emptyset$
and $\overline K\cap \R\subset \{0\}$.
Then Proposition~\ref{p.first} holds with $\Omega^b$ in place of $\Omega$,
and any probability measure on $\Omega^b$ satisfying any
one of conditions 1--4 of the proposition is supported on $\Omega$.
To see this, suppose that $\P$ is a probability measure on $\Omega^b$.
The proofs of the implications
$1\Leftrightarrow 2 \Rightarrow 3 \Leftrightarrow 4$
of Proposition~\ref{p.first}
are valid for $\P$ without modification.
Now suppose that $\P$ satisfies Condition 3 of the proposition.
If $\P[K\text{ connected}]<1$,
then there is a smooth simple path $\gamma:[0,1]\to\overline\H\setminus\{0\}$
such that $\gamma[0,1]\cap\R=\{\gamma(0),\gamma(1)\}$ and
and with positive $\P$-probability
$K\cap\gamma[0,1]=\emptyset$ and $\gamma[0,1]$ separates $K$ in
$\overline\H$.
(This follows because there is a countable collection of 
smooth candidate paths, which is dense, in the appropriate sense.)
Let $A$ be the hull satisfying $\H\cap\p A=\gamma(0,1)$.
If $0\in A$, then $\lim_{s\nearrow 1}\Phi_{\gamma[0,s]}'(0)=0$,
implying $\P[K\cap \p A=\emptyset]=0$, contradicting our
assumptions.  If $0\notin A$, then
$\lim_{s\nearrow 1}\Phi_{\gamma[0,s]}'(0)=\Phi_A'(0)$,
contradicting $\P[K\cap A\ne\emptyset,\ K\cap\p A=\emptyset]>0$.
Hence, $\P[K\text{ connected}]=1$.
Similar arguments show that a.s.\ $0\in\overline K$ and
$K$ is unbounded.  Thus, $\P[\Omega]=1$.
Using this fact, Lemma~\ref{l.uniq} may be applied, giving the
remaining implication $3\Rightarrow 2$.
\end {remark}

\section{Brownian excursions}\label{bexc}

An important
example of a restriction 
measure is given by  
the law of the Brownian excursion from $0$ to infinity in $\H$.
Loosely speaking, this is simply
planar Brownian motion started 
from the origin and conditioned to stay in $\H$ at all 
positive times.
It is closely related to the 
``complete conformal invariance'' of (slightly different)
measures on Brownian excursions in \cite {LW2,LSWup}.

\begin{figure}
\centerline{\includegraphics*[height=2.3in]{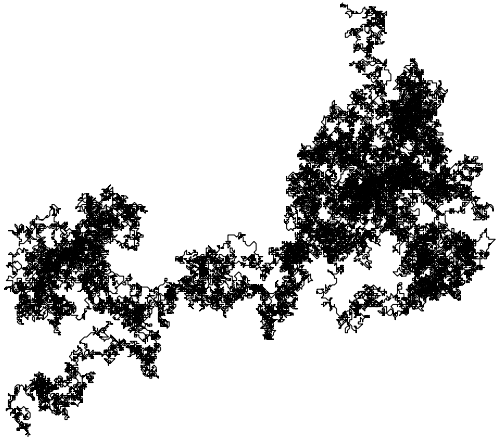}}
\caption{\label{f.exc}A sample of the beginning of
a Brownian excursion.}
\end{figure}
 
Let  
$X$ be a  standard one-dimensional Brownian motion and $Y$ an 
independent three-dimensional Bessel process (see e.g., \cite {RY}
for background  on three-dimensional Bessel processes, its
relation to Brownian motion conditioned to stay positive
and stochastic differential equations). 
Let us briefly recall that a three-dimensional Bessel process is 
the modulus (Euclidean norm)
 of a three-dimensional Brownian motion, and that it can be 
defined as the solution to the stochastic differential equation
$ 
dY_t = dw_t +  dt / Y_t$,
where $w$ is standard Brownian motion in $\R$. It is very easy to see that
$(1/Y_t, t \ge t_0)$ is a local martingale for all $t_0>0$, and 
that if $T_r$ denotes the hitting time of $r$ by $Y$, 
then the law of $(Y_{T_r + t} , t < T_R - T_r )$ is 
identical to that of a Brownian motion started from $r$ and
conditioned to hit $R$ before $0$ (if $0<r<R$).
Note that almost
surely $\lim_{t \to \infty} Y_t = \infty$.
 
The Brownian excursion can be defined as  $B_t = X_t+iY_t$.
In other words,
 $B$ has the same law as the solution to the 
following stochastic differential equation: 
\begin {equation}
\label {e.sde}
dB_t = dW_t + i \frac {1}{ \Im (B_t)} dt
\end {equation}
with $B_0 = z=x + iy$,
where $W$ is a complex-valued Brownian motion,
$x \in \R$ and $ y \geq 0$. 
Note that $B$ is a strong Markov process and that 
$B(0, \infty) \subset \H$ almost surely.
Assume that $y=0$,   and
 let $T_r$ denote the hitting time of the line $\I_r :=
\R + ir$
by  $B$ (i.e., the hitting time of $r$ by $Y$).
Let $S$ denote a random variable with the same 
law as $B_{T_1}$.
Scaling  shows immediately that for all 
$0<r<R$, the law of $B [ T_r, T_R]$ is the law of a Brownian motion
started with the same law as $rS$, stopped at its first hitting
of $\I_R$, and conditioned to stay in $\H$ up to that time.
Note that the probability of this event is $r/R$.

The next proposition, which is due to B\'alint Vir\'ag~\cite{Virag},
implies that the law of the filling $\fillHR(B)$ of the path of
a Brownian excursion $B$ from $0$ in $\H$ is $\P_1$.

\begin{proposition}
\label {p.P1}
 For all  $A \in \hulls^*$,
$\Prob[B[0, \infty) \cap A = \emptyset]  
  = \Phi_A'(0)$.
\end{proposition}

For completeness, a proof is included.

\proof
 Let $\Phi = \Phi_A$.
 Suppose   that  $W$ is a planar Brownian
motion and 
$Z$ is a Brownian excursion in $\H$, both
starting at $z \in \H \setminus A$.
When $\Im (z) \to \infty$, $\Im (\Phi^{-1}(z)) = \Im (z) + o(1)$.  
 Hence, with a large 
probability (when $R$ is large), a Brownian motion started from $z \in \I_R$
(respectively, $z \in \Phi^{-1} (\I_R)$) will hit $\Phi^{-1} (\I_R)$ 
(resp., $\I_R$) before $\R$. 
The strong Markov property  
of planar Brownian motion therefore shows that when $R \to \infty$,
$$  
\P[W \text{ hits } \I_R\text{ before } A\cup\R]  
\sim  
\P[ W \text{ hits } \Phi^{-1}(\I_R)\text{ before } A\cup\R]  
.$$  
But since $\Phi\circ W$ is a time-changed Brownian motion, 
and $\Phi:\H\setminus A\to \H$, 
the right-hand is equal to the probability that a 
Brownian motion started from $\Phi(z)$ hits $\I_R$ before $\R$, namely, 
$ \Im ( \Phi(z)) / R$. Hence, 
$$
\P[Z \text{ hits } \I_R\text{ before } A] 
=
\frac{\P[{W \text{ hits }\I_R\text{ before }A\cup\R}]} 
{\P[W \text{ hits }\I_R\text{ before }\R]} 
=
\frac { \Im [\Phi (z)]}{\Im (z) } + o(1) 
$$
when $R \to \infty$. In the limit $R \to \infty$, we get 
\begin {equation}
\label {e.Im}
\P [ Z  \subset \H \setminus A  ]  
= \frac { \Im [\Phi(z)]}{\Im (z)} = \frac{\Im [g_A(z)]}{\Im(z)}
.\end {equation}
When $z \to 0$, 
$ \Phi (z) = z  \Phi' (0) + O ( |z|^2) $  so that 
$$
\P [ 
 B[0, \infty) \cap   A = \emptyset ] 
= \lim_{s \to 0} \P [ B [ s, \infty ) \subset  \H \setminus A] 
= \lim_{s \to 0} \expect 
\left[ \frac {\Im (\Phi (B_s)) }{ \Im (B_s)} \right] 
= \Phi'(0) 
$$
(one can use dominated convergence here, since $\Im ( \Phi(z)) \le 
\Im (z)$ for all $z$).
\QED

We have just 
proved
that the two-sided restriction measure $\P_1$ 
exists.
By filling unions of $n$ independent excursions, 
one constructs the probability measures $\P_n$ 
for all integers $n \ge 1$ which therefore also exist.
It follows 
(using the fact that the dimension 
of the boundary of the filling of a 
Brownian excursion is $4/3$ \cite {LSWa, Be2}) that for any positive
integer $n$,
the Hausdorff
dimension of the boundary of $K$ defined under
 $\P_n$ is almost surely $4/3$.

\begin{figure}
\centerline{\includegraphics*[height=2.3in]{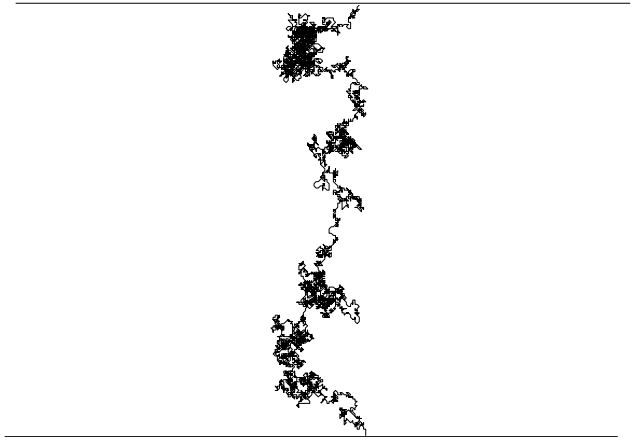}}
\caption{\label{f.excstrip}A sample of a Brownian excursion
from $0$ to $i$ in $\R \times [0,1]$.}
\end{figure}

We have already mentioned that the choice of the 
domain $\H$ and of the boundary points $0$ and
$\infty$ was somewhat arbitrary.
In another simply connected open domain $O \not= \C$
with two invariant boundary points $a$ and $b$ 
(via a given conformal
map $\psi$ from $\H$ onto $O$),
 the Brownian excursion from $a= \psi(0)$ to
$b = \psi (\infty)$ is the solution (up to time-change)
 started from $a$ of 
$$ 
dB_t = d \beta_t + \nabla \varphi (B_t) / \varphi (B_t) dt 
$$
where $\varphi = \Im \,\psi^{-1}$ 
and $\beta$ denotes planar Brownian motion.

\medbreak
Using almost the same proof as in 
Proposition \ref{p.P1} (but keeping track of the  
law of the path),  one can prove the following:
 
\begin{lemma}  \label{mar26.lemma1}
Suppose $A \in \hulls^*$
and $B$ is a Brownian excursion in $\H$
starting at $0$.  Then the conditional law of $(\Phi_A(B(t)), 
t \ge 0)$
given $B \cap A = \emptyset$ is the same as a time
change of $B$.
\end{lemma}

Finally, let us mention the following result that will 
be useful later on.

\begin{lemma}
\label{mar18.prop1}
Let $\Prob^{x+iy}$ denote the law of a Brownian excursion $B$
starting at $x+iy\in\overline{\H}$.
  Then for
every $A \in \hulls^*$,
$$
\lim_{y \rightarrow \infty} y  \,
\Prob^{x+ iy}\bigl[B[0,\infty) \cap A
   \neq \emptyset\bigr] = a(A),
$$
and
\[  \lim_{y \rightarrow \infty}\; y  \int_{-\infty}^\infty
\Prob^{x+ iy}\bigl[B[0,\infty) \cap A
   \neq \emptyset\bigr]\; dx =\pi \; a(A) , \]
where 
$a(A)$ is as in~\eref{adef}.
\end{lemma} 

\noindent
{\bf Proof.}  By (\ref{e.Im}) and the normalization of $g_A$ near 
infinity, 
\begin{eqnarray*}
\Prob^{x+ iy}[B[0,\infty) \cap A
   \neq \emptyset] & = & 1 - \frac  {\Im[g_A(x+iy)]} y \\
& =  & \frac{a(A) \; y}{x^2 + y^2}  + O\left(\frac{1}
   {x^2 + y^2}\right) ,
\end{eqnarray*}
and the lemma readily follows.
\QED

Using Cauchy's Theorem, for example, it is easy to see that the second
statement of the lemma may be strengthened to
\begin{equation}\label{e.cau}
 y  \int_{-\infty}^\infty
\Prob^{x+ iy}\bigl[B[0,\infty) \cap A
   \neq \emptyset\bigr]\; dx =\pi \; a(A),\qquad
y>\sup\{\Im \,z:z\in A\}.
\end{equation}

\section {Conformal image of chordal SLE}
  \label{framesec} 

Let $W:[0,\infty)\to\R$ be continuous with $W_0=0$, 
and let $(g_t)$ be the (chordal) Loewner chain driven by $W$
satisfying~\eref{e.gt}. 
It is easy to verify by differentiation and~\eref{e.gt} 
that the inverse map $f_t(z) = g_t^{-1}(z)$ satisfies 
\[    \p_t f_t(z) = - \frac{2 f_t'(z) }{z - W_t}, \;\;\;\;
                f_0(z) = z . \]

\begin{figure} 
\centerline{
\SetLabels 
(.02*.6)$g_t$\\ 
\B(.45*.95)$G$\\ 
\T(.09*.01)$W_t$\\ 
\L(.97*.6)$\tilde g_t$\\ 
(.5*.3)$h_t$\\ 
\B(.25*.72)$A$\\ 
\B(.3*.06)$A_t$\\ 
\T(.68*.01)$\tilde W_t$\\ 
\endSetLabels 
\AffixLabels{\includegraphics*[width=5in]{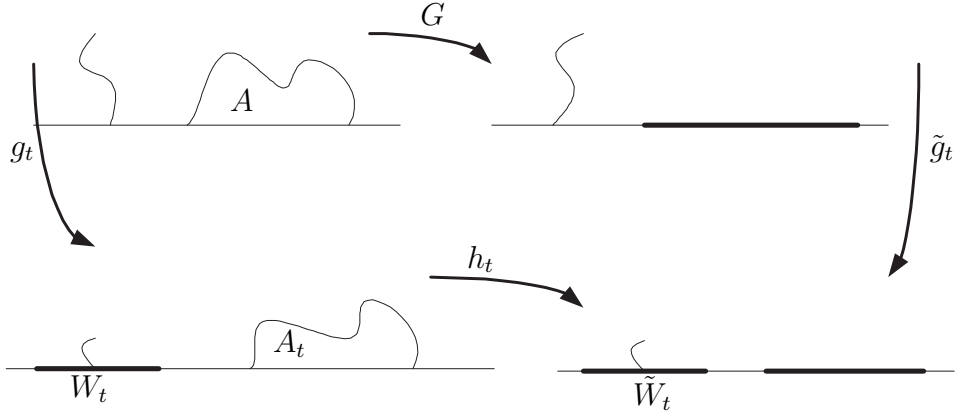}}
} 
\caption{\label{f.comdia}The various maps.} 
\end{figure}

Suppose that  $A \in \hulls^*$ is fixed, and let   $\GG = g_A$  and 
$T = T_A = \inf\{t: K_t\cap A \ne\emptyset \} . $
For $t < T$, let $A_t=g_t(A)$, $\tilde K_t  = \GG (K_t)$ and
$\tilde g_t=g_{\tilde K_t}$.
See Figure~\ref{f.comdia}. 
Then $\tilde g_t$ has an expansion
\[           \tilde g_t(z) = z + \frac{a(t)}{z} + o(z^{-1}) ,
          \;\;\;\; z \rightarrow \infty , \]
where the coefficient $a(t)$ depends on $\GG$ and $W_t$.

  Note that $\tilde g_t$
satisfies the Loewner equation
\[           \p_t \tilde g_t(z) = \frac{\partial_t a(t)}{\tilde g_t(z)
                    - \tilde W_t}, \;\;\;\; \tilde g_0 (z) = z, \]
where 
$$
\tilde W_t := h_t(W_t),\qquad
h_t := \tilde g_t \circ \GG \circ g_t^{-1}=g_{A_t}\,.$$
(This follows from the proof of Loewner's theorem,
because $\tilde g_t(\tilde K_{t+\delta}\setminus \tilde K_t)$
lies in a small neighborhood of $\tilde W_t$ when
$\delta>0$ is small.  Also see \cite[(2.6)]{LSW1}.) 
The identity~\eref{e.aprod} gives
$a(g_t(K_{t+\Delta t}\setminus K_t))
 = 2 \Delta t$.  The image of
$\tilde K_{t+\Delta t}\setminus\tilde K_t$ under $\tilde g_t$
is $h_t\bigl(g_t(K_{t+\Delta t}\setminus K_t)\bigr)$. 
The scaling rule~\eref{ascal} 
of $a$ tells us that as $\Delta t \rightarrow 0+$,
the half-plane capacity of $h_t\bigl(
g_t(K_{t+\Delta t}\setminus K_t)\bigr)$ is
asymptotic to $h_t'(W_t)^2 \cdot 2 \Delta t$.
(The higher order derivatives of $h_t$ can be ignored,
as follows from~\eref{adef2}.  Also see~\cite[(2.7)]{LSW1}.) 
  Hence, \begin{equation}\label{e.aprime}
\partial_t a(t) = 2 h_t'(W_t)^2 .
\end{equation}

Using the chain rule
we get
\begin{equation}  \label{feb14.1}
  [\p_t h_t](z) =  \frac{2 \: h_t'(W_t)^2}{
       h_t(z) - \tilde W_t} - \frac{2\; h_t'(z)}{z - W_t} . 
\end{equation}
This formula is valid for $z \in \H \setminus g_t(A)$
as well as for $z$ in a punctured neighborhood of
$W_t$ in $\R$.  In fact, it is also valid at $W_t$ with
\[  [\p_t h_t] (W_t) = \lim_{z \rightarrow W_t}
     \Bigl( \frac{2 \: h_t'(W_t)^2}{
       h_t(z) - \tilde W_t} - \frac{2 \; h_t'(z)}{z - W_t} \Bigr)
    = - 3 h_t''(W_t) . \]
Computations of a similar nature appear (in a deterministic
setting) in \cite {CM}.
Differentiating (\ref{feb14.1}) with respect to $z$ gives
the equation
\[   [\p_t h_t'](z) = - \frac{2 \: h_t'(W_t)^2\; h_t'(z)}
    {(h_t(z) - \tilde W_t)^2} + \frac{2\; h_t'(z)}{(z - W_t)^2}
                         - \frac{ 2\; h_t''(z)}{z - W_t} . \]
Therefore, at $z = W_t$, 
\[ [\p_t h_t'](W_t)  = \lim_{z \rightarrow W_t}   \p_t h_t'(z)
    =  \frac{h_t''(W_t)^2}{2 h_t'(W_t)} - \frac{4 h_t'''(W_t)}{3} .
    \]
Higher derivatives with respect to $z$ can be handled
similarly.  

Now suppose  that $(B_t, t \ge 0)$ is a standard one-dimensional
 Brownian motion and
that $(W_t, t \ge 0)$ is a
(one-dimensional) semimartingale satisfying $W_0 = 0$
and
\[         dW_t = b_t \; dt + \sqrt \kappa \; dB_t  \]
for some measurable process $b_t$ adapted to the
filtration of $B_t$ which satisfies
$\int_0^t|b_s|\,ds<\infty$ a.s.\ for every $t>0$.

It\^o's formula shows that $\tilde W_t=h_t(W_t)$, $t<T$,
is a semimartingale
with  
\begin{equation}\label{e.dtildeW}
d\tilde W_t = d[h_t(W_t)] =
          h_t'(W_t) \; dW_t + \bigl((\kappa/2) - 3\bigr)  h_t''(W_t)\; dt .
\end{equation}
Here, we need a generalized It\^o's formula since the function
$h_t$ is random (see, e.g., exercise (IV.3.12) in \cite {RY}).
However, since $h_t'(z)$ is $C^1$ in $t$, no
extra terms appear.  Similarly,
\[ d[h_t'(W_t) ] = h_t''(W_t) \; dW_t + \Bigl(
         \frac{h_t''(W_t)^2}{2 h_t'(W_t)}  + \bigl(\frac \kappa 2 
     - \frac 4 3 \bigr) h_t'''(W_t)\Bigr) \; dt .
\]
Let $\aa > 0$ and let $Y^0_t = h_t'(W_t)^\aa$.  Then yet another application
of It\^o's formula gives   
\begin{equation}
\label{dYt}
 \frac{dY^0_t}{\aa \, Y^0_t} =    \frac{h_t''(W_t)}{h_t'(W_t)} \; dW_t  
   +  \Bigl( \frac{(\aa-1) \kappa + 1}{2} \: \frac{h_t''(W_t)^2}
    {h_t'(W_t)^2}   + (\frac\kappa 2 - \frac 4 3 )
     \frac{h_t'''(W_t)} {h_t'(W_t)} \Bigr) \; dt . 
\end{equation}
These computations imply readily the following results:

\begin{proposition}\label{p.loc}
Let $b_t=0$ (i.e., $W_t = \sqrt \kappa B_{t }$)
and fix $A\in\hulls^*$.
\begin {enumerate}
\item
$\tilde W_{t}$
is a local martingale
if and only if  $\kappa = 6$.
\item
Suppose that $\kappa=6$
and let $T':=\inf\bigl\{t:K_t\cap \Phi_A(\partial A)\ne\emptyset\bigr\}$.
Then $(\Phi_A(K_t),\; t<T)$, has the same law as a time change
of $(K_t,\; t<T')$.
\end {enumerate}
\end{proposition}

Claim 2 is basically the \lq\lq locality property\rq\rq, which is central
to the papers~\cite{LSW1,LSW2,LSW3}.  It has been proven in~\cite{LSW1}
using a somewhat different, longer and more technical, proof.
(See~\cite{LSW1} for a more complete discussion of this important property.)
Further consequences of this locality result are discussed in 
\S \ref {s.fileq}.

\proof
Statement 1 is clear from~\eref{e.dtildeW}.
To prove 2, set $\kappa=6$ and
define $\tau(s):=\inf\{t:a(t)\ge 2\,s\}$,
$\hat W_s:=\tilde W_{\tau(s)}$ and
$\hat g_s:=\tilde g_{\tau(s)}$.
Then by~\eref{e.aprime}
$\partial_s \tau(s)=2/ \partial_{\tau} a(\tau)=h_t'(W_t)^{-2}$.
Hence,~\eref{e.dtildeW} shows that
$\hat W_s/\sqrt 6$, $s<\tau^{-1}(T)$, is stopped Brownian motion.
Statement 2 now follows since
$$
\p_s \hat g_s(z)=\frac 2{\hat g_s(z)-\hat W_s}\,.
$$
\QED

\begin{proposition}
\label {p.rest}Suppose $b_t = 0$ (i.e., $W_t = \sqrt \kappa B_t$),
$\aa > 0$, $\kappa>0$.
The process $Y^0_t = h_t'(W_{t})^\aa$, $t<T$, is a local  martingale 
 for all $A \in \hulls^*$ if and
only if $\kappa = 8/3$ and $\aa=5/8$.
\end{proposition}
\proof Immediate from~\eref{dYt}.
\QED

The next section will be devoted to consequences of this property
of SLE$_{8/3}$.
 
Before stating a useful generalization of Proposition~\ref{p.rest},
we recall a few basic facts about
the Schwarzian derivative,
\[   Sf(z) := \frac{f'''(z)}{f'(z)} - \frac{3 f''(z)^2}{2f'(z)^2} . \]
An essential property of the Schwarzian derivative is
that $Sm=0$ when $m$ is a M\"obius transformation,
$m(z)=(az+b)/(cz+d)$, $ad-bc\ne 0$.
An easy direct calculation shows that
 $Sf(0)= -6\, a(A)$ when $f(z)=g_A(-1/z)$.
Since $S(m\circ f)=S(f)$ for M\"obius transformations $m$,
it follows that 
\begin{equation}\label{e.aschw}
S g_A(0) = -6\,a\bigl( \{-z^{-1}:z\in A\}\bigr).
\end{equation}
Consequently, $Sg_A\le 0$ on $\R\setminus A$.
(In fact, in \S\ref{bubsec}
we show that $-Sg_A(z)/6$ is  a hitting measure for
Brownian bubbles.)
If $\alpha >0, \lambda \in \R$, let
\[ \Yl_t :=   h_t'(W_t)^\alpha 
 \; \exp\Bigl(\lambda \int_0^t \frac {Sh_s(W_s)}{6} \; ds\Bigr). \]

\begin{proposition} \label{mar19.prop1}
If $W_t = \sqrt \kappa B_t$, 
\begin {equation}
\label {defalpha}
\alpha  = \frac {6 - \kappa}{2\kappa}
\end {equation}
 and
\begin {equation}
\label {deflambda}
\lambda = \frac {(8 - 3 \kappa)(6 - \kappa)}{2 \kappa},
\end {equation}
 then $\Yl_t$, $t<T$, 
is a local martingale.  If $\kappa \leq 8/3$, then $\Yl_t$ is a bounded 
martingale (in fact, $0 \leq \Yl_t \leq 1$). 
\end{proposition}

\noindent
{\bf Proof.} The local martingale property follows immediately
from It\^o's formula and (\ref {dYt}). The bound on $\Yl_t$ follows 
from $Sh_t(W_t) \leq 0$ and~\eref{e.positive}.
\QED

\section{Restriction property for SLE$_{8/3}$}
\label {s.83}
We now discuss some  consequences of Proposition \ref {p.rest}
for SLE$_{8/3}$. In particular, we establish the following
theorem.
 
\begin{theorem}[Restriction]  \label{restheorem}
Let $\gamma$ be the SLE$_{8/3}$ path starting
at the origin and $A \in \hulls^*$, then
\[    \Prob [\gamma[0,\infty) \cap A = \emptyset] 
            = \Phi_A'(0)^{5/8} . \]
The law of $\gamma (0, \infty)$ is therefore $\P_{5/8}$.
\end{theorem}

Roughly, we will need to show that as $t\nearrow T$
the $Y^0_t$ of Proposition~\ref{p.rest} 
converges to $0$ or $1$, respectively, if
$\gamma$ hits $A$ or not.
Some simple deterministic 
lemmas below will help us establish this.

\begin{lemma}\label{gotoone} 
Let $A\in\hulls_+$, let $W:[0,\infty)\to\R$ be continuous, 
and let $g_t$ be the corresponding solution of~\eref{e.gt}. 
Let $K_t$ be the associated growing hull, as defined 
in Section~\ref{s.prelim}, 
and suppose that $\bigcup_{t>0} K_t\cap A=\emptyset$.
Let
$T(r) : = \sup\bigl \{ t \ge 0  :  K_t\subset r\U\bigr  \}$
and $A_t=g_t(A)$.
Then
$$\lim_{r \to \infty}
g_{A_{T(r)}}'(W_{T(r)})=
 1\,.
$$
\end{lemma}

\proof
Set $T:=T (r)$, $W:=W_T$,
$a_0:=\inf (A\cap\R)$, $a_1:=\sup(A\cap\R)$, $A':=A\cup[a_0,a_1]$
and $\tilde A:= g_T(A')$.
Let $r'=\sup\{|z|:z\in K_T\}$ (actually, $r'=r$, but we don't
need this fact), and let $z_0\in K_T$ be such that $|z_0|=r'$.
Set $\beta(s):=z_0+(1-s)i$, $s\in[0,1]$.
Then the limit $w:=\lim_{s\nearrow 1}g_T\circ\beta(s)$
exists.  (This is because the image of the conformal map $g_T$ is
a smooth domain, i.e., $\H$.  See, e.g., \cite[Proposition 2.14]{PommBB}.)
Moreover, since $\beta(1)=z_0$, we must have $w=W$.
Otherwise, one easily gets a contradiction to
$z_0\in K_T\setminus\bigcup_{t<T} K_t$.
The extremal length (see~\cite{A2} for the definition
and basic properties of extremal length) from $A'$ to
the circle $|z|=r'\ge r$ goes to infinity with
$r$.  By monotonicity and conformal invariance of
extremal length, the extremal length in $\H$ between
$\tilde A$ and $(-\infty,W]$ goes to infinity as well.
Because $\tilde A$ is connected, this implies that
$\diam(\tilde A)/\inf\{|W-z|:z\in \tilde A\}$
goes to zero.  Since $g_B'(W)$ is invariant under
scaling $B$ about $W$, under translating $B$ and $W$,
and is monotone decreasing in
$B$, this means that when $r$ is large,
$g_{\tilde A}'(W)$ is at least $g_B'(0)$,
where $B=\{z\in\overline\H:|z-1|\le\eps\}$ and
$\eps>0$ is arbitrarily small.  The lemma follows.
\QED

\begin{lemma}\label{deter}
Let $W:[0,\infty)\to\R$ be continuous,
let $g_t$ be the corresponding solution of~\eref{e.gt}, 
and let $K_t$ be the associated growing hull, as defined 
in Section~\ref{s.prelim}. 
Let $A\in\hulls^*$ be a smooth hull.
Suppose that
$T:=\inf\{t\ge 0:K_t\cap A\ne\emptyset\}<\infty$
and $K_T\cap A\cap\R=\emptyset$.
Set $A_t:= g_{t}(A)$, $t<T$.
Then
$$
\lim_{t\nearrow T} \Phi_{A_t}'(W_t)=0\,.
$$
\end{lemma}
\proof 
We first argue the rather obvious fact $K_T\cap A\subset\p A$.
Note that $\lim_{t\searrow 0}\diam(K_t)=0$ follows from 
the continuity of $W_t$ at $t=0$.  This implies that for
every $z\in\H$ the harmonic measure of $K_t$ from $z$ in
$\H$ goes to $0$ when $t\searrow 0$.
For $s>0$, the hull evolution
$t\mapsto \overline{ g_s(K_{t+s}\setminus K_s)}$
is driven by $t\mapsto W_{s+t}$.
By conformal invariance of harmonic measure,
since $W_t$ is uniformly continuous in $[0,T]$,
for every $z\notin K_T$ the harmonic measure
of $K_{t+s}\setminus K_s$ in $\H\setminus K_s$
goes to zero uniformly as $t\searrow 0$ while $s\in[0,T]$.
In particular, the harmonic measure
of $K_T\setminus\bigcup_{t<T} K_t$ in
$\H\setminus\overline {\bigcup_{t<T} K_t}$ is zero
from any $z$ in the latter open set.
Since $K_T$ does not contain $\p A$, we may apply
this from a point in $\p A\setminus K_T$,
and conclude that $K_T\cap A\subset\p A$.

Let $z_0$ be some point in $\p A\cap K_T$.
Let $\beta:[0,1]\to\H$ be a smooth path
such that $\beta[0,1)$ is contained in the interior
of $A$, $\beta(1)=z_0$, and $\beta'(1)$ is
orthogonal to $\partial A$ at $z_0$.
By smoothness of $\p A\cap\H$, there is some small disk $D\subset\H$
with center $z_0$ such that $\p D\cap\p A$ consists of exactly two points.
Let $\sigma_1$ be the arc of $\p A$ from $z_0$ to
$\p D$ that goes away from $z_0$ in the direction
$i\beta'(1)$, and let $\sigma_2$ be the other arc of $\p A$
from $z_0$ to $\p D$.

We claim that 
$$
\hat\beta(x):=g_{T}\circ \beta(x)-W_T\,,\qquad x\in[0,1)\,,
$$
is a path which is contained in a sector
$\bigl|\Re\, z\bigr|\le c\,\Im\, z$ for some $c$.
As in the proof of Lemma~\ref{gotoone}, we know
that $\lim_{x\nearrow 1} \hat\beta(x)=W_T$.
Suppose that we start a two dimensional Brownian motion
from $\beta(x)$, $x\in[0,1)$, and stop when we hit $K_T\cup\R$.
Then there is probability bounded from below
that we hit this set to the ``right'' of $\beta[0,1]$,
because the Brownian motion has probability bounded from 
below to first hit $\p D\cup \sigma_1$ on $\sigma_1$
and from the side of $\sigma_1$ not in the interior of $A$.
(The careful reader might want to draw a little figure here.)
By conformal invariance, this shows that 
Brownian motion started from $\hat\beta(x)$ has probability bounded from
below to first hit $\R$ in $(0,\infty)$.
Consequently, $ \Re\, \hat\beta (x)\ge - c\,\Im\, \hat\beta(x)$
for some constant $c$ independent of $x$.  The symmetric argument
also shows $ \Re\,\hat\beta (x)\le c\,\Im\, \hat\beta(x)$, for some $c$,
and the claim is established. 

Since
$$
\hat\beta(x)=\lim_{t\nearrow T}g_{t}\circ \beta(x)-W_t\,,
\qquad x\in[0,1)\,,
$$
it follows that when $t$ tends to $T$ a Brownian motion excursion from
$0$ to $\infty$ in $\H$ will hit $\hat\beta$ before
exiting $g_{t}(D\setminus K_t)-W_t$ with probability tending to $1$.
This implies that this Brownian motion excursion will hit
$A_t-W_t$ with probability tending to $1$. 
By Proposition~\ref{p.P1},
this probability is the same
as $1-\Phi_{A_t}'(W_t)$.
\QED

\proofof{Theorem \ref{restheorem}}
By Proposition~\ref{p.first}, it suffices to consider the case where $A$
is a smooth hull in $\hulls_+\cup\hulls_-$.
By symmetry, we may take $A\in\hulls_+$.
Proposition~\ref{p.rest} shows that $Y^0_t = 
h_t'(W_{t })^{5/8}$ is a bounded continuous local
martingale.  By the martingale convergence
theorem, the a.s.\ limit $Y^0_T := \lim_{t \nearrow T} Y^0_t$ 
exists and $Y^0_0 = \expect [Y^0_T]$, 
where $T=\sup\{t:\gamma[0,t]\cap A=\emptyset\}$.
Lemmas~\ref{gotoone} and \ref{deter} show that
$Y^0_T=1_{T=\infty}$ a.s. 
This proves the theorem.
\QED

Combining this with the results of Section \ref {rest}
shows the following identity in law:

\begin {corollary}\label{c.5o8}
The filling of the union of 8 independent chordal SLE$_{8/3}$'s 
has the same law as the filling of the union of 5
independent Brownian 
excursions from $0$ to infinity.  In both cases, the law
is $\P_5$, the two-sided restriction measure with
exponent $5$.
\end {corollary}

Theorem \ref{restheorem} suggests that SLE$_{8/3}$ should be the 
limit as the lattice mesh goes to zero
of the self-avoiding walk.
See \cite {LSWsaw} for a discussion of these conjectures.
Also, we know \cite {RS} that SLE$_{8/3}$ is a simple
curve. This suggests that $\alpha = 5/8$ is the smallest 
possible value for which $\P_\alpha$ exists. We shall see
later (Corollary \ref{mar26.cor1}) that this is indeed the case.

\medbreak

A similar proof to that of Theorem \ref{restheorem},
using Proposition \ref{mar19.prop1}, gives the following
important generalization.

\begin{theorem}  \label{mar19.theorem1}
Suppose $0 \leq \kappa \leq 8/3$
and let $\alpha$ and $\lambda$ be as in (\ref {defalpha})
and (\ref {deflambda}).
If\/ $W_t = \sqrt{\kappa} B_t$ and $A \in \hulls^*$, then
\begin{equation}\label{e.withbub}
    \Phi_A'(0)^{\alpha}  = \E \;  \left[ \; 
     1_{\{\gamma[0,\infty) \cap A = \emptyset \}} \;
    \exp\Bigl( \lambda \int_0^\infty \frac {Sh_s(W_s)}{6} \; ds \Bigr)  \; 
   \right] . \end{equation}
\end{theorem}

The following section will provide a more concrete meaning to the
right hand side and use the theorem to construct the measures
$\P_\alpha$, $\alpha \ge 5/8$.

\section {Bubbles}
\label {s.bub}
\subsection {Brownian bubbles}  \label{bubsec}

We now define the measure $\nu$ of
{\em Brownian bubbles hanging at infinity}.
 This is a $\sigma$-finite but infinite measure
on unbounded closed connected sets $K\subset\H$ such that
$\C\setminus K$ is connected.
The definition of $\nu$ is rather simple.
For $z\in\overline\H$, let $\P^z$ denote the law of the
Brownian excursion $Z_t$ started at $z$, as discussed in
\S\ref{bexc}.  Let $\hat \P^z$ denote law of the filling of
$Z$, $\fillHR(Z[0,\infty))$.
Set
\begin{equation}\label{nudef}
\nu :=\frac 1\pi \lim_{y\to\infty} y\, \int_{\R} \hat \P^{x+iy}\,dx\,.
\end{equation}
In other words, one considers the limit when $y \to \infty$ of the (infinite) measure 
obtained by filling a Brownian excursion that is started on the line $\I_y$, where the initial point is chosen according to $y/\pi$ times the Lebesgue measure on $\I_y$.
The existence of the limit is easily justified, as follows.
Set $\sigma_y = \inf \{ t   :\Im\, Z_t=y \}$.
Recall that for $y'>y$, $\P^{x+iy'}[\sigma_y<\infty]=y/y'$.
On the event $\sigma_y<\infty$, let $Z^y$ be the path
$Z^y_t = Z_{\sigma_y+t}$, $t\ge 0$.
By the strong Markov property it therefore follows that
the image of the measure
$1_{\sigma_y<\infty}\int \P^{x+iy'}\,dx$
under the map $Z_t\mapsto Z^y_t$ is precisely
$(y/y')\int \P^{x+iy}\,dx$.
The existence of the limit~\eref{nudef} readily follows.

Suppose that $A \in \hulls$.
We have by Lemma~\ref{mar18.prop1}
\begin {equation}
\label {e.in}
\nu [ K \cap  A \not= \emptyset ] = 
\frac 1\pi\lim_{y\to\infty} y \int_{\R}
\P^{x + i y} [ Z \cap  A 
\not= \emptyset  ]  \; dx 
  =a(A) .
\end {equation}
This can be used to 
give an alternative proof of the existence of the limit in
the definition of $\nu$.

Let $A$ and $A'$ in $\hulls$. Define
$A'' = A \cup g_A^{-1} (A')$. Then $A'' \in \hulls$ and 
$a( A'') = a( A) + a(A')$ by~\eref{e.aprod} and invariance of
$a(A)$ under real translations.
Hence, 
\begin {eqnarray*}
 {\nu [ g_A (K) \cap A' \not= \emptyset
\hbox { and } K \cap A = \emptyset  ] }
&=&
\nu [ K \cap A'' \not= \emptyset  \hbox { and }
K \cap A = \emptyset ]
\\
&=&
\nu [ K \cap A'' \not= \emptyset ] - \nu [ K \cap A \not= \emptyset ]
\\
&=&
a ( A'') - a (A) = a(A').
\end {eqnarray*}
Therefore, the image of $1_{K\cap A=\emptyset}\,\nu$
under $g_A$ is $\nu$.  In the terminology of \S\ref {s.prelim},
this says that $\nu$ is invariant under the semigroup 
$\{g_A:A\in\hulls\}$.

\medbreak

Now define the measure $\mu$ on Brownian bubbles at $0$
as the image of $\nu$ under the  inversion
$z \mapsto -1/z$. It is a measure on the set of bounded $K\subset\H$
with $\overline K=K\cup\{0\}$.
By~\eref{e.aschw} we have for $A\in\hulls^*$
\begin{equation}\label{muschw}
\mu[K\cap A\ne\emptyset]=-\frac {Sg_A(0)}6\,.
\end{equation}
We may think of $\mu$ as a measure on the space
$\Omega_b$ of connected bounded sets $K\subset\H$
such that $\overline K=K\cup\{0\}$ and $\C\setminus \overline K$ is
connected.

If $\Gamma$ denotes the semigroup of dilations $r_\lambda$,
$r_\lambda(z) = \lambda z$,
 then it is easy to see from~\eref{ascal} that $\nu$ and
$\mu$ are 
$\Gamma$-covariant; in fact, $r_\lambda \circ \nu =
\lambda^{-2} \nu, r_\lambda \circ \mu = \lambda^{2} \nu$.
 
\medbreak

There are a number of alternative equivalent ways to define 
the measures $\nu$ and  $\mu$ and to derive their properties:
\begin {itemize}
\item
Define a measure on two-sided excursion in $\H$, $Z_t$, starting
at the origin by 
\[   Z_t = \left\{ \begin{array}{ll} Z^1_{-t}, & - \infty < t \leq 0, \\
                       Z^2_t,& 0 \leq t < \infty, \end{array} \right. \]
where $Z^1,Z^2$ are independent excursions in $\H$
starting at $0$.
The measure $\nu$ is obtained by choosing a point
$z \in \H$ according to two-dimensional
Lebesgue measure and letting
$B_t = Z_t + z$. 
(This definition gives $\nu$ as a measure on paths, rather than
fillings.  The parametrization of the paths is chosen so that
at time $0$ they attain their minimal imaginary part.)
\item
Define
on different
spaces a one-dimensional Brownian excursion
$(e_t, 0\le t \le T)$ (defined under
It\^o's excursion measure $n$) and a
Brownian
bridge $( b_t,0\le t \le 1)$ with $b_0= b_1 = 0$ (defined
under a probability measure $P$).
Recall that It\^o's excursion measure is obtained as the
limit when $\eps \to 0$ of $\eps^{-1}$ times the
probability measure defining a Brownian motion started
from $\eps$ and killed at its first hitting time of 0, $T$.
We then define the process
$$
Z (t) = T^{1/2}\, b (t/T) + i\, e_t,\qquad t \in [0, T].
$$
If $\Upsilon$ denotes the map $(b,e)\mapsto \fillHR (Z)$,
then $\mu =
c\, \Upsilon( T^{-1/2} n \otimes P)$ for some constant $c$.
The factor $T^{-1/2}$ is needed in order for $\mu$ to scale
properly under the dilations $r_\lambda$.

\item
We can also relate Brownian bubbles to Brownian excursions in $\H$.
Given $\eps>0$, let $m_\eps:\H\to\H$ denote an arbitrary
M\"obius transformation fixing $0$ and satisfying $m_\eps(\infty)=\eps$.
Then $\mu$ is (the filling of) the limit when $\eps \to 0$
of $\eps^{-2}$ times
the $m_\eps$-image of the law of a Brownian excursion from $0$
to $\infty$ in $\H$.

\end {itemize}

\begin {remark}
The previous description can be combined with the fact that
eight SLE$_{8/3}$ are equivalent to five Brownian excursions
(i.e., more precisely, Corollary~\ref{c.5o8})
to describe the measure on Brownian bubbles using
SLE$_{8/3}$. Actually, since we only focus on the
hull of the Brownian bubbles, the description of its
frontier in terms of SLE$_{8/3}$ is 
natural.
The probability that a SLE$_{8/3}$ (or a Brownian excursion)
in $\H$ from $0$ to $\eps$ hits the circle of radius
$\delta$ around zero decays like a constant times
$\eps^2$ when $\eps$ goes to zero.
Hence, if we condition the union of eight SLE$_{8/3}$
(resp., five Brownian excursions) 
from $0$ to $\eps$
to intersect this circle, and take the limit when $\eps \to 0$,
we obtain exactly the same outer boundary (we know from
the previous description that this limit
exists) as if we condition just one SLE$_{8/3}$ or one Brownian
excursion, since with high probability only
one of them will hit the circle. Hence, we get that $\mu$ is
the filling of
the limit when $\eps \to 0+$ of $8\eps^{-2} /5$ times the law of chordal
SLE$_{8/3}$ in $\H$ from $0$ to $\eps$.
Hence, the hull of a Brownian bubble is also an
``SLE$_{8/3}$-bubble''.
\end {remark}

\begin {remark}
Let $t>0$ and
let $X_t$ be a sample from the Poisson point process with mean (intensity)
$t\,\nu$.
\footnote{\label{f.pois}This means that $X$ is a countable random set of
bubbles such that whenever $D_1,D_2,\dots,D_k$ are disjoint measurable
sets of bubbles the random variables
$|X\cap D_j|$, $j=1,\dots,k$, are independent
and $\E\bigl[|X\cap D_j|\bigr]=\nu[D_j]$.}
Let $U_t$ be the filling of the union of bubbles in $X_t$,
$ U_t=\fillHR \bigl(\bigcup X_t\bigr) $.
By the properties of Poisson point processes, for all $A\in\hulls$
$$
\P [ U_t  \cap A = \emptyset ] 
= \exp \bigl(- t\, \nu [ \{ K  :  K \cap A \ne \emptyset \} ]\bigr)
=
\exp\bigl (- t \,a(A)\bigr). 
$$
But for $A,A'\in\hulls$ we have $a(A\cdot A')=a(A)+a(A')$,
by~\eref{e.aprod}.  Hence, the law of $U_t$ is
covariant under the semigroup $\{\Phi_A:A\in\hulls\}$.
However, it is not scale-invariant, because the image of 
$\nu$ under $z \mapsto \lambda z$ is $\lambda^{-2} \nu$.
Thus, the distribution of $U_t$ under the map $z \mapsto \lambda
z$ is the same as that of $U_{\lambda^{-2} t}$.
This shows that the assumption
of $\Gamma$-invariance in Statement~1 of Proposition~\ref{p.first}
is important.
\end {remark}

\subsection {Adding  a Poisson cloud of bubbles to SLE}

We are now ready to give a rather concrete interpretation of
the right hand side of~\eref{e.withbub} and thereby construct the
measures $\P_\alpha$ for $\alpha\ge 5/8$.

Suppose that $\kappa \le 8/3$ and, as in \S\ref{framesec},
let
\[ \alpha = \alpha_\kappa = \frac{6 - \kappa}{2\kappa}, \;\;\;\;
\lambda = \lambda_\kappa = \frac{(8 - 3 \kappa)
(6-\kappa)}{2 \kappa}\,. \]
Consider a Poisson point process $X$ on
$\Omega_b\times [0,\infty)$ with mean (intensity)
$\lambda\,\mu\times dt$,
where $dt$ is Lebesgue measure (see footnote~${}^{\ref{f.pois}}$
for the definition).
As before, let $\gamma$ denote the SLE$_\kappa$ path,
$g_t$ the corresponding conformal maps,
and $W_t$ the Loewner driving process.
We take
$\gamma$ to be independent from $X$.
Since $\kappa\le 8/3$, we know 
from \cite {RS}  that $\gamma$ is a simple curve.
Let 
$$
\hat X:=\bigl\{ g_t^{-1}(K+W_t): (K,t)\in X,\, t\in [0,\infty)\bigr\},
$$
and let $\Xi$ be the filling of the union of elements of $\hat X$
and $\gamma$,
$$
\Xi=\Xi({\kappa}) := \fillHR\Bigl(\gamma(0,\infty)\cup\bigcup\hat X\Bigr).
$$

Let $A\in\hulls^*$, and let $h_t$ be the normalized conformal map from
$\H \setminus {g_t(A)}$ onto $\H$ as in \S\ref{framesec}.
By~\eref{muschw}, for any $t>0$ on the event $\gamma[0,t]\cap A=\emptyset$,
\begin {eqnarray*}
\lefteqn {\P\bigl[\{K:g_t^{-1}(K+W_t)\cap A\ne\emptyset\}\bigm|g_t\bigr]
}
\\
&=&
\P\bigl[\{K:(K+W_t)\cap g_t(A)\ne\emptyset\}\bigm|g_t\bigr] 
\\
&=&
-Sh_t(W_t)/6\,,
\end {eqnarray*}
where $K$ is independent from $\gamma$ and has law $\mu$.
Consequently, on the event $\gamma[0,\infty)\cap A=\emptyset$,
$$
\P\bigl[\Xi\cap A=\emptyset\bigm|\gamma\bigr]=
 \exp \Bigl(
  \lambda \int_0^{\infty} \frac {Sh_t (W_t)}{6}  \, dt\Bigr)\,.
$$
By taking expectation and applying Theorem \ref{mar19.theorem1},
we get
\begin{equation}\label{xilaw}
\P\bigl[\Xi\cap A=\emptyset\bigr]=
\Phi_A'(0)^\alpha\,,
\end{equation}
which almost proves,

\begin {theorem}
\label {loops}
For any $\kappa \in[0, 8/3]$, the law of\/ $\Xi(\kappa)$ 
 is $\P_{\alpha_\kappa}$.
\end {theorem}
\proof
Given the discussion above, all that remains is to
show that $\overline \Xi=\Xi\cup\{0\}$.
Let $D:=\{z\in\overline\H:|z-x_0|\le\eps\}$, where $0<\eps<1$ and
$x_0\in[1,2]$.
Then $1-\Phi_D'(0)=O(\eps^2)$.  Consequently, 
$\P\bigl[\dist(\Xi,[1,2])<\eps\bigr]=O(\eps)$.  Thus, a.s.,
$\overline\Xi\cap[1,2]=\emptyset$.  By scaling, it follows that
$\overline\Xi\cap(\R\setminus\{0\})=\emptyset$ a.s.
Let $X_{t_0}^{t_1}$ denote the set of pairs $(K,t)\in X$
with $t\in [t_0,t_1)$.
Since the $\mu$-measure of the set of bubbles of diameter larger
than $\eps$ is finite, a.s., for every $t_1\in[0,\infty)$ the set
of $(K,t)\in X_{0}^{t_1}$ such that $K$ has diameter at least $\eps$ is finite.
Therefore, the set
$\gamma[0,t^1]\cup\bigcup\bigl\{ g_t^{-1}(K+W_t): (K,t)\in X_0^{t_1}\bigr\}$
is closed a.s.\ when $t_1<\infty$.
To show that $\Xi\cup\{0\}$ is closed, it therefore suffices to prove
that $\bigcap _{s>0}\overline \Xi_s=\emptyset$, where
$\Xi_s:=\gamma(s,\infty)\cup\bigcup
\bigl\{ g_t^{-1}(K+W_t): (K,t)\in X_s^{\infty}\bigr\}
$.
Let $T(R)$ denote the first time $t$ such that $|\gamma(t)|=R$, 
and let $A=\{z\in\overline \H:|z|\le 1\}$.  Let $A_+$ denote the 
set of points in $A$ which are to the right of $\gamma[0,T(R)]$ 
or on $\gamma$
(i.e., the intersection of $A$ with the
closure of the domain bounded by $\gamma[0,T(R)]\cup[0,R]$ and  
an arc of the semicircle $\{z\in\H:|z|=R\}$).
The proof of Lemma~\ref{gotoone} gives $\Phi_{g_{T(R)}(A_+)}'(W_{T(R)})\to 1$ 
as $R\to\infty$.
The stationarity property of SLE with
equation~\eref{xilaw} imply
$$
 \P\bigl[\Xi_{T(R)}\cap A_+=\emptyset\bigm| \gamma[0,T(R)]\bigr] 
=
\Phi_{g_{T(R)}(A_+)}'(W_{T(R)})^\alpha\to 1\,. 
$$
A symmetric argument shows that this holds with $A$ in
place of $A_+$.  This implies that
a.s.\ $\bigcap _{s>0}\overline \Xi_s$ is disjoint
from the disk $|z|<1$.
Scale invariance now gives
$\bigcap _{s>0}\overline \Xi_s=\emptyset$ a.s., and completes the proof.
\QED

The theorem shows that for all $\alpha > 5/8$, the
measure $\P_\alpha$ exists and can be constructed by adding
bubbles with appropriate intensity to SLE$_{\kappa}$
with $\kappa = 6 / (2 \alpha +1 )$.
The frontier of the set defined under $\P_\alpha$ 
has Hausdorff
dimension $4/3$ (because of the Brownian bubbles).
For instance, for integer $\alpha$, this shows that
SLE$_\kappa$ can be coupled with
the union of $n$ independent excursions so as to be a subset of their filling.

Note that $\lim_{\kappa\to 0+}\lambda_\kappa=\infty$,
while $\lambda_{8/3}=0$.
Also observe that $\alpha_2=1$, so that
adding Brownian bubbles to SLE$_2$ with appropriate density
gives the measure on hulls of Brownian excursions.
This is not surprising
since  SLE$_2$ is
the scaling limit of loop-erased simple
random walks as proved in \cite {LSWlesl}.

In~\cite {LSWloops} it is shown that there
is a natural Poisson point process $L$ of sets in $\H$,
independent from $\gamma$, such
that $\Xi$ can be also described as the (filling of)
the union of $\gamma$ with  those sets in $L$ which meet $\gamma$.

\section {One-sided restriction}
\label {s.os}

\subsection {Framework}

 Recall the definition of $\hulls_+$
from \S\ref{s.prelim}.
Set ${\cal A}_+=\{\Phi_A:A\in\hulls_+\}$ .
Let $\Omega_+$ denote the set of 
all closed connected sets $K \subset \overline \H$ 
such that $K \cap \R = (-\infty, 0]$ and
$\H \setminus K$ is connected.
We endow $\Omega_+$ with the $\sigma$-field generated by the 
family of events 
$\{ K \cap A = \emptyset\}$, where $A \in \hulls_+$.
We say that the probability measure $\P$ on $\Omega_+$
satisfies the {\em right-sided restriction 
property} if it is $\mathcal A_+$ covariant
 and scale invariant.
In other words,
$\P[K\cap(A\cdot A')=\emptyset]=
\P[K\cap A=\emptyset]\,\P[K\cap A'=\emptyset]
$
and
$\P[K\cap A=\emptyset]=\P[K\cap (\lambda A)=\emptyset]$
hold for all $A,A'\in\hulls_+$, $\lambda>0$.

The proof of Proposition~\ref{p.first} shows that
if $\P$ satisfies the right-sided restriction
property, then there exists a constant $\alpha\ge 0$
such that for all $A \in\hulls_+$,
$$ \P [ K \cap A = \emptyset] = \Phi'_A(0)^{\alpha}.$$
Conversely,
for all $\alpha\ge 0$, there exists at
 most one such probability measure $\P$. If
it exists, we call it the
{\em right-sided restriction measure with exponent $\alpha$} and
denote it by $\P_\alpha^+$.  For $\alpha \geq 5/8$,
we may obtain $\P_\alpha^+$ by applying $\fillneg$
to a sample from the
two sided restriction measure $\P_\alpha$. 
(Recall the notation $\fillneg$ from \S~\ref{s.prelim}.)

In the following, we will see two other constructions
of $\P_\alpha^+$, which are valid for all $\alpha>0$,
the first is based on reflected Brownian motion,
while the second is an SLE type construction, where
an appropriate drift is added to the driving process of SLE$_{8/3}$.
We will also be able to conclude that $\P_\alpha$ does
not exist when $\alpha<5/8$.

We generally ignore the uninteresting case
$\alpha=0$, where $K=(-\infty,0]$ a.s.

\subsection {Excursions of reflected Brownian motions}

We now construct
$\P_\alpha^+$ for all $ \alpha > 0$ using reflected Brownian
motions, or,  more precisely,   Brownian
excursions conditioned to avoid $(0,\infty)$ and reflected
at angle $\theta$ off
$(-\infty,0]$. 
In order to define this,
fix $\theta \in (0, \pi)$ and let $c = c_\theta = -\cot \theta$.
We first consider
Brownian excursions in the wedge  
$$W:= W(\theta) = \{ r e^{i \varphi} \ : \ r > 0 \hbox { and } \varphi \in (0, 
\pi - \theta) \}.$$
reflected in the horizontal direction off the boundary line $x = c\,y$.
Let $(Y_t, t \ge 0)$ denote a three-dimensional Bessel process started from $0$
(i.e., a one-dimensional Brownian motion conditioned to
stay in $(0,\infty)$).
Let $(X_t, t \ge 0)$ be a one-dimensional Brownian motion,
acting independently of $Y$  when inside  $W(\theta)$, 
which is reflected off   the ``moving boundary'' $c Y_t $. More 
precisely, if $\tilde X$ is standard Brownian motion started from $0$, 
independent of $Y$,
then $X_t$ is the unique continuous function such that $X_t \ge c Y_t$ and
$X_t  = \tilde X_t + \ell_t$, where $\ell$ is a non-decreasing continuous function 
with $\int 1_{X_t> c Y_t } d\ell_t = 0$.
(See, e.g., \cite {RY} for more on Skorokhod's reflection lemma).

Define
$$
Z_t = X_t + i Y_t.
$$
Let $V$ denote a random variable which has the same 
law as
$Z_{T_1}$
where $T_R$ denotes for all $R>0$ the hitting time of $R$ by $Y$.
Then, for all $r<R$, the process $X+iY$ on $[T_r,T_R]$ is started 
with the same  distribution as $rV$ and then evolves like 
two-dimensional Brownian motion which is reflected horizontally
off the line $y=c\,x$, and conditioned to hit $\R +i\,R$ before $\R$ (this 
event is independent of $X_{T_r}$ and has probability $r/R$ for
the unconditioned reflected Brownian motion). 

Let $b(z) = b_\theta(z) = z^{\pi/(\pi - \theta)}$.
If  $A \in \hulls_+$, let  
$  F = F_{A,\theta} = b^{-1} \circ \Phi_{b(A)} \circ b $. Then $F$ is a 
  conformal
transformation of $W \setminus A$ onto $W$ with
$F(0) = 0$ and $|F(z) - z|$ bounded.
It is straightforward to show
that the image under $F$ of a 
horizontally reflected Brownian motion in $W$, up to the first time
it hits   $\R  \cup  (W \setminus A$),  is a (time-changed)
horizontally reflected Brownian motion in $W$:
as long as it is away from the line $y=c\,x$, this is just 
conformal invariance of planar Brownian motion, and 
since $F'$ is real on the line $y=cx$, it follows
that $F(W)$ also gets a horizontal push when it hits the line;
that is,  $d\tilde \ell_t = F'(Z_t) \,d \ell_t$  defines
the corresponding push for $\Re\bigl(F(Z_t)\bigr)$.
It therefore follows just as in the case of the Brownian excursion in $\H$
(which   corresponds
to the limiting case $\theta=0$) 
that for all small $r>0$, 
$$
\P \bigl[ Z[T_r, \infty) \cap  A = \emptyset\bigr]
=
\frac {\E [ \Im ( F (Z_{T_r}) )]}{r} \,.
$$
Hence, letting $r \to 0$, we get by dominated convergence that 
$$
\P [ Z \cap  A  = \emptyset ] = \Phi_{b(A)}'(0)^{1 - (\theta/\pi)}
\,.
$$

We now define the ``reflected Brownian excursion'' in $\H$ (in short 
RBE) as $B=\{b_\theta(Z_t):t\ge 0\}$.
Then, the previous equation for $\P[Z\cap A=\emptyset]$ may be rewritten
$$
\P [
B \cap  A = \emptyset  ] =
\P [ Z \cap  b_\theta^{-1}(A) = \emptyset ]  
=
\Phi_A'(0)^{1- (\theta/ \pi)},
$$
which shows that $\fillneg(B)$
satisfies right-sided restriction with exponent $\alpha=1- \theta /\pi$.

Note that the limiting cases $\theta=0$ and $\theta= \pi$ correspond
respectively to the Brownian excursion ($\alpha=1$) and to the 
ray $(-\infty, 0]$ that stays on the boundary ($\alpha=0$). 

\medbreak

Reflected Brownian excursions therefore show that for all $\alpha \in (0,1]$,
the right-sided restriction measure with exponent $\alpha$ exists. 
Taking unions of independent hulls which satisfy the 
right-sided restriction property, yields a realization of 
another right-sided 
restriction measure (and the exponent add up).  We summarize
this in a proposition.

\begin{proposition}
\label {rexc}
  The right-sided restriction measures
$P_\alpha^+$ exist for all $\alpha > 0$.  If 
$\alpha = a_1 + \cdots + a_k$ where $k$ is a positive integer and $a_1,
\ldots, a_k 
\in (0,1]$, then $\fillneg$ applied to the union of
$k$ independent RBEs with respective angles
$\theta_1 = \pi(1-a_1), \ldots , \theta_k
= \pi (1 - a_k)$ has law 
$P_\alpha^+$.
\end{proposition}

This, together with the observation that $P_\alpha^+$
can be realized as the ``left-filling''
of samples from $P_\alpha$ (when the latter exists)
   implies various rather surprising identities in law between 
``right-boundaries'' of different processes:

\begin{corollary}$\;$
\begin {enumerate}
\item
The right-boundary of an RBE   
with angle $3 \pi / 8$ 
has the same law as   SLE$_{8/3}$.
In particular, its law is symmetric with respect to the 
imaginary axis.

\item
The right boundary of the union of $n$ independent RBE with angles 
$\pi - \theta_1, \ldots, \pi - \theta_n$
 has the same law as the right boundary of an RBE with 
angle $\pi - (\theta_1+\cdots + \theta_n)$,
provided that $\theta_1+ \cdots + \theta_n < \pi$.
 
\item
The right boundary of the union of two independent RBE which are
orthogonally reflected on the negative half-axis 
has the same law as the right boundary of a Brownian excursion.
\end {enumerate}
\end{corollary} 

The first statement
shows that the Brownian 
frontier (outer boundary) looks like a locally symmetric path.
This, in spirit, answers a question raised
 by Chris Burdzy after Benoit Mandelbrot
noted (based on simulations) the similarity between the dimension 
of self-avoiding walks and the Brownian frontier and
proposed \cite {Ma}
the name ``self-avoiding Brownian motion'' for the 
Brownian frontier.
Burdzy's question was whether the Brownian frontier is
[locally] symmetric.  There are several different precise formulations
of this question.  See Section~\ref{rem} for more about this issue.

\begin{figure}
\centerline{\includegraphics*[height=2.3in]{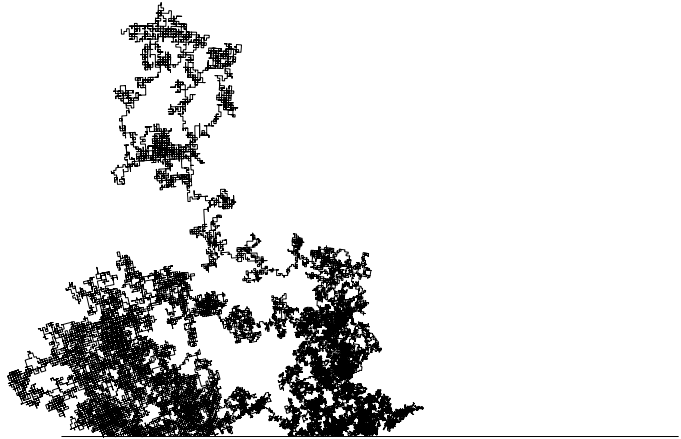}}
\caption{\label{f.rbe}A sample of the beginning of 
a (long) reflected random walk
 with angle $3\pi/8$ (its right boundary looks like SLE$_{8/3}$).} 
\end{figure}
 
Note that the last two statements (and their proofs) do not  
use SLE.
The first statement yields an extremely fast algorithm to simulate 
chordal SLE$_{8/3}$ and therefore also the scaling limit
of self-avoiding half-plane walks (modulo the
conjecture \cite {LSWsaw} that chordal SLE$_{8/3}$ is the scaling
limit of the half-plane self-avoiding walk) as the right-boundary of 
a reflected excursion. See \cite {Kennalgo} for an
 algorithm 
to simulate directly  such walks.
 
\subsection {The SLE(${\kappa, \rho}$) process}
\label{slekapparho}

We will now describe  the 
right-boundaries of these sets  
 in terms of SLE-type
paths that are driven by Bessel-type processes.

Before introducing these processes, let us give 
a brief heuristic.
Let $\gamma$ denote the right-boundary of a 
Brownian excursion in $\H$.
Let us condition on a piece $\gamma[0,t]$.
For the future of $\gamma$ beyond time $t$,
the right-hand boundary of $\gamma[0,t]$ acts just
like the positive real axis, $\gamma[t,\infty)$ is ``conditioned'' not
to hit the right hand side of $\gamma[0,t]$. 
If we believe in conformal
invariance of the process, then we may ignore all
the geometry of the domain $\H\setminus\gamma[0,t]$
and map it onto the upper half plane.  However, we
should keep track of the left image of $0$ under
the uniformizing map $g_t$.  It is reasonable to believe that
this is all that would be relevant to the distribution
of $g_t(\gamma[t,\infty))$. 
(We will a posteriori see that this is the case.)
  Let $W_t=g_t(\gamma(t))$, let
$O_t$ be the left image of $0$ under $g_t$, and
take $t$ to be the half-plane capacity parametrization for $\gamma[0,t]$.
Then the pair $(W_t,O_t)$ is a continuous Markov process,
and the chordal version of Loewner's theorem gives
  $dO_t/dt=2/(O_t-W_t)$.
  Scale  and translation invariance show that it is enough to
know what happens to $W$ infinitesimally when
$O=0$ and $W=1$.  The natural guess is that at that moment 
we have $dW=\sqrt{\kappa}\,dB+\rho \,dt$, for some constants 
$\kappa>0$ and $\rho\in\R$.  Scaling this to other values of $W$ gives 
SLE($\kappa, \rho$), as will be defined shortly.

Suppose that $ \kappa > 0$, $\rho > -2$ and that $B_t$ is a standard
one-dimensional Brownian motion.
Let $(O_t,W_t)$ be the solution of 
\[    dO_t =  \frac{2 \; dt}{O_t - W_t} , \;\;\;\;
            dW_t = \frac{\rho \; dt}{W_t - O_t} + \sqrt{\kappa} \; d B_t \]
with $O_0 = W_0 = 0$ and $O_t \le W_t$.
The meaning of this evolution is straightforward at times
when $W_t>O_t$, but a bit more delicate when $W_t=O_t$.
  One way to construct $(O_t,W_t)$ is
to first define $Z_t$ (later to become $W_t-O_t$)
as the solution to the Bessel equation
\[   dZ_t = \frac{(\rho +2) \; dt}{Z_t} + \sqrt \kappa \; dB_t  \]
started from $Z_0 = 0$.
More precisely,  $Z_t$ is $\sqrt \kappa$ times a $d$-dimensional
Bessel process where 
\begin{equation}\label{ddef}
d= 1 + \frac {2(\rho+2)}{\kappa}\,.
\end{equation}
It is well-known (e.g., \cite {RY})
 that this process is well-defined (for all $\rho
> -2$ and all $t \ge 0$).
 Note also that $\int_0^t du / Z_u = (Z_t - \sqrt {\kappa} B_t ) / (
\rho +2 ) < \infty$ for all $t \ge 0$.
Then, set
\begin {eqnarray*}
    O_t &=&  - 2 \int_0^t \frac{du}{Z_u}\,, \\ 
     W_t &=&  Z_t + O_t \,. 
\end {eqnarray*}
If we then define the family of conformal maps $g_t$ by 
$\partial_t g_t (z) = 2 ( W_t - g_t(z))^{-1}$ and $g_0 (z) = z$
(for $z \in \overline \H$),
we get a Loewner chain that we call chordal SLE($\kappa, \rho$).
Note that  when $\rho= 0$, we get the ordinary chordal SLE$_\kappa$.
Intuitively, the definition of SLE ($ \kappa, \rho$) can be understood as
follows: $O_t$ is the left-most point
of $g_t(\partial K_t)$ (when $K_t$
  is a simple path, this is simply the ``left'' image of the origin under 
$g_t$); the $W_t$ gets a push away from this point
if $\rho >0$ (or towards this point if $-2<\rho < 0$),
 and this push is ``constant'' modulo scaling.

The next lemma lists a few basic properties of SLE($\kappa,\rho$),
which are generalizations of known results for SLE; i.e., for the
case $\rho=0$.

\begin {lemma} 
\label {l.propsle}
Let $\kappa>0$, $\rho>-2$ and set $\rho_0:= -2 + \kappa/2$.
Let $K_t$ denote the evolving hulls of
SLE($\kappa,\rho$) and $K_\infty:=\bigcup_{t>0} K_t$.
\begin {enumerate}
\item
The distribution of 
SLE($\kappa, \rho$) is scale-invariant.
More precisely, if $\lambda>0$, then
$(K_t,t\ge 0)$ has the same distribution as
$(\lambda^{-1}\,K_{\lambda^2 \,t},t\ge 0)$.
\item
If $\kappa\le 4$ and $\rho\ge\rho_0$, then a.s.\ 
$K_\infty \cap\R=\{0\}$.
\item
If $\kappa\le 4$ and $\rho <\rho_0 $, then a.s.\ 
$K_\infty \cap\R=(-\infty,0]$.
\item $K_\infty$ is a.s.\ unbounded.
\end {enumerate}
\end {lemma}

Recall that a.s.\ the $d$-dimensional
Bessel process returns to zero if and only if $d<2$.
This will be essential in the proof of 2 and 3.

\proof
Clearly, $(W_t, O_t)_{t \ge 0}$ has the same scaling
property as Brownian motion, and 1 follows.

Now assume $\kappa\le 4$.
Let $\tau_1=\sup\{t\ge 0:1\notin K_t\}$.
We want to show that $\tau_1=\infty$ a.s.
Set $x_t=g_t(1)$ for $t< \tau_1$ and
observe that $x_t-O_t$ is monotone increasing.
In particular, $x_t-O_t\ge 1$, $t<\tau_1$.
On the set of times $t<\tau_1$ such that
$x_t-W_t<1/2$, we therefore have
$d W_t \le 2\, |\rho|\,dt+\sqrt\kappa\,dB_t$.
Setting $\tilde x_t=x_t-W_t$, we get
$$
d\tilde x_t \ge - 2\,|\rho|\,dt-\sqrt\kappa\,dB_t+ (2/\tilde x_t)\,dt\,,
$$
on the set of times $t$ such that $\tilde x_t<1/2$.
If $\rho= 0$, by comparing with the Bessel process we see
that a.s.\ $\tilde x_t$ never hits $0$ and so $\tau_1=\infty$.
For $\rho\ne 0$, note that for any finite fixed $t_0>0$ and any
$c\in\R$ the law of the process $(B_t + c\,t, t\le t_0)$
is equicontinuous with the law of $(B_t, t\le t_0)$.
(In fact, after conditioning on the position of the process
at time $t_0$, their distribution is identical.)
Therefore, also in this case $\tilde x_t$ never hits $0$
and $\tau_1=\infty$.  Hence, a.s.\ $1\notin K_t$
for all $t\ge 0$.  This also implies that
$K_t\cap [1,\infty)=\emptyset$ a.s.\ for all $t\ge 0$,
since $K_t\cap\R$ is an interval.
Scale invariance then implies
$K_\infty \cap (0,\infty)=\emptyset$ a.s.

Now suppose $\rho\ge \rho_0$.
Then the Bessel process $Z_t/\sqrt\kappa$
has dimension $d\ge 2$, as given by~\eref{ddef}.
Consequently, a.s.\ $W_t-O_t=Z_t>0$ for all $t>0$.
If $x<0$, then $g_t(x)\le O_t$ for all $t\ge 0$.
Hence, $K_\infty\cap (-\infty,0)=\emptyset$ a.s.

Now take $\rho\in(-2,\rho_0)$.
Set $y_t = g_t (-1)$ for
$t< \tau_{-1}:=\sup\{t\ge 0:-1\notin K_t\}$.
Using $W_t - y_t \ge W_t - O_t$ and $\rho < 0$, we get for $t<\tau_{-1}$
\begin {eqnarray*}
W_t - y_t 
&=&
\sqrt {\kappa} B_t + \int_0^t ds \left( \frac {\rho}{W_s-O_s} + 
\frac {2}{W_s-y_s}
\right)
\\ & \le &
\sqrt {\kappa} B_t + \int_0^t \frac {(\rho +2) \,ds}{W_s - y_s}\,.
\end {eqnarray*}
So that $W_t - y_t$ is smaller than a Bessel process that hits 
zero a.s. Hence, a.s.\ $-1\in K_\infty$.
This implies $[-1,0]\subset K_\infty$ a.s.,
and by scaling $(-\infty,0]\subset K_\infty$ a.s.
This completes the proof of 2 and 3.

Statement 4 easily follows from 1, for example.
One could also use the fact that the half-plane capacity of
$K_t$ is $2t$.
\QED

The SLE($8/3, \rho$)'s are related to the measures $\P_\alpha^+$
via the following theorem that will be proved in the 
next subsection.

\begin {theorem} 
\label {onesided}
Let $\rho>-2$, and let $K=\fillneg\bigl(\overline K_\infty\bigr)$,
where $K_t$ is the hull of SLE($8/3,\rho$) and
$K_\infty=\bigcup_{t\ge 0} K_t$.
Then $K$ satisfies the right-sided restriction
property with exponent
$$
\alpha = \frac {20 + 16 \rho + 3 \rho^2}{32} 
= \frac {(3 \rho +10)(2 + \rho) }{32}  
. $$
\end {theorem}

Note that when $\rho$ spans $(-2, \infty)$, $\alpha$ spans $(0, \infty)$.
This theorem has several nice corollaries, some of which we now briefly 
discuss.

\begin{corollary}
If $\alpha \geq 5/8$, the
 right boundary of the two-sided restriction measure
$\P_\alpha$ has the same law as the SLE($ 8/3 ,
\rho (\alpha)$) path, where
$$\rho (\alpha) =  \frac {-8 + 2 \sqrt { 24 \alpha +1 } } 3 .$$
In particular, the right boundary of a Brownian excursion has
the law of (the path of) SLE($8/3,2/3$) and the right boundary of the union
of two Brownian excursions has the law of SLE($8/3,2$).
\end{corollary}

\begin {corollary}  \label{mar26.cor1}
For all
 $\alpha<5/8$, the two-sided restriction probability measure $\P_\alpha$
does not exist.
\end {corollary}

\proof
Note that when $\rho<0$, $W_t - \sqrt {\kappa} B_t$ is decreasing. It 
follows easily that the probability that $i$ ends up 
eventually to ``the right'' of the right hand boundary of SLE($8/3 , \rho$)
(i.e., $i$ is separated from $1$ by
$K_\infty\cup(-\infty,0]$),
 is strictly larger than the corresponding quantity for SLE($8/3,0$),
which is $1/2$ by symmetry.
However, the same symmetry argument 
shows that for any $\alpha>0$, if the two-sided 
probability measure with exponent $\alpha>0$ exists, then the 
$\P_\alpha$ probability that $i$ ends up to the ``right''
of $K$ is at most $1/2$ (it can be smaller if $K$ is of positive 
Lebesgue measure). If $K$ has law $\P_\alpha$ with some
$\alpha < 5/8$, then $\fillneg(K)$ has law
$\P_{\alpha}^+$, which is described
using SLE ($ 8/3, \rho$) for  some $\rho <0$. This contradicts the 
fact that the probability that it passes to the left of $i$ 
is at least $1/2$.
\QED

\begin{corollary}\label{realcutpts}
 The boundary of the right-sided 
restriction measure intersects the negative half-line if and only 
if $\alpha< 1/3$. In particular, the reflected Brownian excursion 
with reflection angle $\theta$ on the negative half-line has 
cut points on the negative half-line if and only if 
$\theta > 2 \pi /3$.  
\end{corollary}

\noindent
{\bf Proof.}
This is just a combination of Lemma \ref {l.propsle}, Theorem \ref {onesided}
and Proposition \ref {rexc}. 
\QED

\begin {remark}\label{r.nocut}
Note that non-existence of cut points on the negative 
half-line for the angle $2 \pi /3$ proves (via the 
correspondence between reflected Brownian motion and the 
SLE$_6$ hull 
that is discussed
in \S \ref {s.fileq}) 
non-existence of cut-points for the  SLE$_6$ hull on the 
positive and negative half-line (and therefore also non-existence
of double points for SLE$_6$ that are also local cut-points for the SLE$_6$ 
path). 
In the discrete case (i.e., critical site 
percolation on the triangular grid), van den Berg and Jarai \cite {BJ} 
have recently derived a stronger version of this result (with decay rates
for probabilities). 
\end {remark}

\subsection{Proof of Theorem \ref{onesided}}
 
Fix $\rho > -2$ and let 
$$c = \frac {3 \rho}8 \hbox { and }  b = \frac {\rho (4 + 3 \rho)}{32}
.$$
Let $(O_t,W_t)$ generate an
  SLE($8/3,\rho$) process 
 so that
\[   dW_t = \frac{\rho}{W_t - O_t} \; dt + \sqrt{8/3} \; dB_t , \;\;\;\;
\  dO_t = \frac{2}{O_t - W_t} \; dt .\]
Let $A \in \hulls_+$ be a given smooth hull,
 and let $\Phi = \Phi_A$, $T= T_A$,
and $h_t$ be as in  \S\ref{framesec}
and define (for $t<T$),
\[  M_t  :=    h_t'(W_t)^{5/8} \; h_t'(O_t)^b \;
     \left[ \; \frac{h_t(W_t) - h_t(O_t)}{W_t - O_t} \; \right]^{c} . \]
Of course, when $W_t=O_t$, we take $M_t=h_t'(W_t)^{(5/8)+b+c}$.

\begin{lemma}  
$(M_t, t <T)$ is a local martingale.
\end{lemma}

\noindent
{\bf Proof.}
Plugging the explicit choice of $(W_t, t \ge 0)$
into the results of 
\S \ref{framesec}  
shows that for $t<T$, 
\begin {eqnarray*}
  d[h_t(W_t)] &=& 
   \Bigl( \frac{\rho\; h_t'(W_t) }{W_t - O_t} 
    -(5/3)  \; h_t''(W_t)\; \Bigr) \, dt  + \sqrt{8/3}\;
 h_t'(W_t) \; dB_t, \\ 
 d[h_t'(W_t) ] &=&  \Bigl(
   \frac{\rho\; h_t''(W_t) }{W_t - O_t} + 
         \frac{h_t''(W_t)^2}{2 h_t'(W_t)}  \Bigr) \, dt  
     + \sqrt{8/3}\;
 h_t''(W_t) \; dB_t, \\ 
   d[h_t(O_t)] 
&=&
 \frac{2 h_t'(W_t)^2}{h_t(O_t) -
                h_t(W_t)}  \,dt , \\ 
   d[h_t'(O_t)] 
&=&
 \Bigl(
 \frac{2 \; h_t'(O_t)}{(O_t - W_t)^2}
-\frac{2 \; h_t'(W_t)^2 \; h_t'(O_t)}{(h_t(O_t) - h_t(W_t))^2} 
\Bigr) \, dt \,. 
\end {eqnarray*}
Using these expressions in It\^o's formula for $dM_t$, one
can now compute the semi-martingale decomposition of $M_t$.
  This is tedious but 
straightforward, so we omit the detailed calculation here.
 The drift term of $dM_t$ turns
out to be $M_t$ times
\begin {eqnarray*}
&&(\frac {5\rho}8 - \frac {5c}3 ) \frac { h_t'' (W_t)}{h_t' (W_t) (W_t -O_t)}
\\&+& 
(2b - c (\rho+2) + \frac 43 c (c+1) ) \frac {1}{(W_t -O_t)^2} \\
&+&
(-2b + 2c + \frac 43 c (c-1) ) \frac { h_t'(W_t)^2 }{( h_t (W_t)- h_t (O_t))^2}
\\
&+&
(-\frac 83 c^2 + \rho c) \frac { h_t' (W_t)}{ (W_t- O_t)(h_t (W_t) - 
h_t (O_t))}.
\end {eqnarray*}
The terms in $h_t''(W_t)^2 / h_t'(W_t)^2$ 
and in $h_t'' (W_t) / (h_t (W_t)- h_t (O_t))$ happen to  
vanish because of the choice of the exponent $5/8$ (and $\kappa= 8/3$).
The lemma follows as this drift term vanishes for the appropriate choice
of $b$ and $c$.
\QED

\begin {lemma}
\label {l.decay}
There exists $\eps>0$ such that 
$M_t \le h_t'(W_t)^\eps$ for all $t<T$. In particular, $M_t \le 1$.
\end {lemma}

\noindent 
{\bf Proof.}
 When $\rho \ge 0$,
the statement is trivial since $b , c \ge 0$ and
 $h_t'(W_t)$, $h_t'(O_t)$ and 
$(h_t (W_t) - h_t (O_t))/(W_t - O_t)$ are all in $[0,1]$. 
One has to be a little bit careful when $\rho <0$ as $c<0$ and 
$b$ can be negative as well.
Let
\[  \alpha = \frac{5}{8} + b + c   = \frac{(3 \rho + 10)(2 + \rho)}{32} \]
and note that $\alpha>0$.

We now want to show that
\begin {equation}
\label {ineq2}
h_t' (W_t) \le \frac { h_t (W_t) - h_t (O_t)}{W_t- O_t} \le h_t'(O_t) \le 1\,.
\end {equation}
This will be established by showing that $h_t'$ is 
decreasing in $(-\infty,W_t]$.  Recall that $h_t=g_{g_t(A)}$.
In fact, the following argument shows that $g_{A^*}'$ is monotone
decreasing on $x < \inf (A^*\cap\R)$ for every smooth hull $A^*$.
Applying this with $A^*=g_t(A)$ then yields~\eref{ineq2}.
To prove this monotonicity result, we realize the map
$g_{A^*}=g_{g_t(A)}=h_t$ as a map in a Loewner chain, as follows.
Let $\beta:[0,S]\to\overline\H$ be the smooth path
$\overline{\partial A^*\cap\H}$
starting from $\beta(0)=\inf(A^*\cap\R)$
and parametrized by half-plane capacity from $\infty$.
Set $\hat g_s=g_{\beta[0,s]}$ and $x_s=\hat g_s\bigl(\beta(s)\bigr)$.
By the chordal version of Loewner's theorem,
$\partial_s \hat g_s (z) = 2/ (\hat g_s (z) - x_s)$.
Then $\hat g_S = h_t $, since both are equal to the normalized map
from $\H \setminus g_t (A)$ onto $\H$.
Since
$
\partial_s \hat g_s'(z) = {-2 \hat g_s'(z)}/{(\hat g_s(z)- x_s)^2}
$, it follows that
\begin {equation}
\label {e.deri}
\partial_s \log \hat g_s' (z) = \frac {-2} { (\hat g_s(z) - x_s)^2}
\,.
\end {equation} 
Therefore, $h_t'$ is decreasing on $x<\beta(0)$,
which proves~\eref{ineq2}.
This, 
implies that 
in the case where $  \rho \in [ -4/3, 0)$ (because
then $b\le 0$ and  $c\le 0$), 
$$
M_t \le h_t'(W_t)^{5/8 + b +c } = h_t'(W_t)^\alpha \le 1.
$$

Now suppose $\rho \in (-2, -4/3)$, which gives $c<0$ and $b>0$. 
For this case, we use a similar argument
involving the Loewner chain $(\hat g_s)$, but a little more care
is necessary.
Suppose that $o<w<x_0$ and let 
$w_s = \hat g_s (w)$, $o_s = \hat g_s (o)$.
{}From the expression for $\p_s \hat g_s$ we have
$$
\partial_s \log ( w_s - o_s)
= \frac {-2}{(x_s - w_s ) (x_s -o_s )}
.$$
Combining this with (\ref {e.deri}) shows that 
\begin {equation}
\label {e.M}
 M_t
=
 \exp\int_0^S 
\left(
 \frac {-2\,(5/8)}{(x_s - w_s)^2}
 + \frac {-2\,c}{(x_s-w_s)(x_s - o_s)}
 + \frac {-2\,b}{(x_s- o_s)^2 }  
\right)ds\,,
 \end {equation}
where $w=W_t$ and $o=O_t$.
But $o_s < w_s  < x_s$ for all $s \le S$. Hence 
(using the explicit values of $b$ and $c$,
and $\rho \in (-2, -4/3)$), 
$$
\frac {({5}/{8}) - \eps}{(x_s  - w_s)^2}
 + \frac {c}{(x_s-w_s)(x_s - o_s)}
 + \frac {b}{(x_s- o_s )^2 }  
\ge 0
$$
for some positive $\eps=\eps (\rho)$.
The lemma finally follows.
\QED

\medb
{\bf End of the proof of Theorem \ref {onesided}.}
It now remains to study the behavior of the bounded 
martingale $(M_t, t<T)$ when $t \to T-$.
Let $K_t$ be the growing hull of the SLE($8/3,\rho$) process.
Recall that $A$ is a smooth hull in $\hulls_+$
and that $K_\infty\cap(0,\infty)=\emptyset$.
Hence, if $T<\infty$, then $K_T\cap  A\ne\emptyset$,
 Lemma \ref {deter} shows that 
$$
\lim_{t \to T-} h_t' (W_t) = 0,
$$
and  Lemma~\ref {l.decay} 
implies that $\lim_{t\to T-}M_t = 0$ if $T< \infty$.
Let us now suppose that 
$T= \infty$. 
Lemma~\ref {gotoone}
shows that a.s.\ on the event $T=\infty$,
$$
\lim_{r \to \infty}  
h_{T_r}' (W_{T_r}) = 1.
$$
By (\ref {ineq2}), it follows 
that $\lim_{r\to \infty} M_{T_r} =1$.
Hence, since $M_t$ converges a.s. and in $L^1$ when $t \to T$,
it follows that 
the limit is $1_{T=\infty}$ and
$$
\P [K_\infty \cap A = \emptyset ] 
= M_0 = \Phi_A'(0)^\alpha.
$$

It remains to prove that a.s.\ $K_\infty\cap A\ne\emptyset$
if and only if  $\overline K_\infty \cap A\ne\emptyset$.
As $K_t$ is closed for each $t$, the proof of this fact is
essentially identical to the argument showing that
$\bigcap_{s>0}\overline\Xi_s=\emptyset$
given at the end of the proof of Theorem~\ref{loops}.
\QED

\subsection{Formal calculations}

   In this subsection
we discuss  how one can
guess the form of the 
martingales $M_t$ giving the intersection probabilities.
Since this is not part of the proof, we will not be
rigorous; however, much of 
 this discussion can be made rigorous 
and may be used to further study restriction measures.

Let $\hulls_1$ denote the set of $A\in\hulls$ such that
$A\cap\R\subset (1,\infty)$.
For $A\in\hulls_1$, 
let $\phi_A$ denote the unique conformal map $\phi:\H\setminus A\to\H$
which fixes each of the three points $0,1,\infty$.
Suppose that $X$ is a random set in $\H$, whose law is
covariant with respect to the semigroup $\Lambda=\{\phi_A:A\in\hulls_1\}$.
An example of such a set should be given by an SLE($8/3, \rho$) started 
with $O_0=0$ and $W_0=1$.

One can also associate to $A$ the unique conformal map 
$g_A :\H \setminus A
\to \H$ that is normalized at infinity.
Note that 
$$ 
\phi_A (z) = (g_A (z) - g_A(0)) / (g_A (1) - g_A(0)).
$$
Define now
$$
\tilde H (g_A) := H ( \phi_A)  
:=\P[X\subset \phi^{-1}(\H)]\,.
$$
Our goal is to show that 
$\tilde H ( g_A) $ is of the form 
$$g_A'(0)^a g_A'(1)^b (g_A(1)-g_A(0))^c.$$
It is more convenient to work first with $H$ since $\Lambda$ is a semi-group
while the family $\{ g_A : A \in \hulls_1\}$ is not. 
The function $H$ is a semigroup homomorphism from $\Lambda$
into the multiplicative semigroup $[0,1]$.  Consequently,
$dH$ is a Lie algebra homomorphism into $\R$.
The \lq\lq basic\rq\rq\ vector fields generating $\Lambda$
have the form
$$
A(x)= \frac{z(1-z)}{z-x}\,, \qquad x>1\,.
$$
(This vector field corresponds to an infinitesimal slit at $x$.
Note that flowing along $A(x)$ preserves $0,1,\infty$.)
This is a one real-parameter ($x$) family of vector fields in the $z$-plane.
The commutator of $A(x)$ and $A(y)$ turns out to be
$$
[A(x),A(y)]= A(x) \p_z A(y)- A(y)\p_z A(x)
= \frac{(x-y)\,(z-1)^2\,z^2}{(x-z)^2\,(y-z)^2}
\,.
$$
This is supposed to be annihilated by $dH$, since $[0,1]$
is commutative.  Hence, if we divide by $x-y$ and take
a limit as $y\to x$, it will also be annihilated
by $dH$.  This is the vector field
\begin{equation}\label{hA}
\hat A(x):=
\lim_{y\to x} (x-y)^{-1} [A(x),A(y)]=
\frac{(1-z)^2\,z^2}{(x-z)^4}\,.
\end{equation}
To understand $H$, we want to determine the function
$$
h(x)=dH\bigl(A(x)\bigr)\,.
$$
So we want to extract from $dH\bigl(\hat A(x)\bigr)=0$ information
about $dH\bigl(A(x)\bigr)$.  For this, we write
$\hat A(x)$ as a linear combination of the derivatives
$\p_x^j A(x)$ with coefficients functions of
$x$.  Direct computation gives:
$$
\hat A(x) =
-\p_x A(x)+ \frac 12\,(1-2\,x)\p_x^2 A(x)
+\frac 16\,(x-x^2)\,\p_x^3 A(x)\,.
$$
Since $dH$ is linear, it commutes with $\p_x$,
and we get
\begin{align*}
0&
=dH\bigl(\hat A(x)\bigr)
\\&
=
-\p_x dH\bigl( A(x)\bigr)+ \frac 12\,(1-2\,x)\,\p_x^2  dH\bigl(A(x)\bigr)
+\frac 16\,(x-x^2)\,\p_x^3  dH\bigl(A(x)\bigr)
\\&
=
-h'(x)+ \frac 12\,(1-2\,x)h''(x)
+\frac 16\,(x-x^2)\,h'''(x)
\,.
\end{align*}
The general solution of this equation turns out to be very simple,
it is
$$
h(x)= \frac {c_0+c_1 \,x+ c_2\,x^2}{x\,(1-x)}
= a_0 \left(\frac { x-1}{x}\right) + a_1 + a_2 \left(\frac {x}{x-1}\right).
$$
This, in fact, already determines the general form of $H$, since any
$\phi_A$ can be obtained in a Loewner-equation way from
the infinitesimal fields $A$.

We now want to translate this information in terms of $\tilde H$,
since this is the framework that we are working with (even though the
present analysis shows that it is not the most natural one
here, but 
we 
have some formulas worked out already, so it is more
economical at this point).
Suppose that $g_A$ is obtained via a Loewner chain driven by a
continuous function $(x_s, s \le S)$, then 
$\partial_s \hat g_s = {2}/ {(\hat g_s (z) - x_s)}$ and $\hat g_S=g_A$.
Associate to each $\hat g_s$ the corresponding function
$\phi_s=(\hat g_s-\hat g_s(0))/(\hat g_s(1)-\hat g_s(0))$,
which is normalized at $0,1,\infty$.
Then,  
$$
\partial_s \phi_s=
\frac {-2}{\bigl(\hat g_s(0)- x_s\bigr)\bigl(\hat g_s(1)-x_s\bigr)}\,
 A\Bigl(\frac{x_s - \hat g_s(0)}{\hat g_s(1)- \hat g_s(0)}\Bigr) \circ 
\phi_s.
$$
Since $H$ is multiplicative,
$\p_s H(\phi_s)=H(\phi_s)\,dH\bigl((\p_s\phi_s)\circ\phi_s^{-1}\bigr)$.
It therefore follows readily that 
$\tilde H(\hat g_S)=H(\phi_S)=\int_0^S \p_s H(\phi_s)\,ds$
is equal to 
$$
\exp \int_0^S
ds \left( \frac {a'}{(x_s -\hat g_s(1))^2} + 
\frac {c'}{(x_s - \hat g_s(0))(x_s - \hat g_s
(1))}
+ \frac {b'}{(x_s -\hat g_s(0))^2 }  \right),
$$
as we had in~\eref{e.M}.
The ``good'' values of $a', b', c'$ can then 
be determined by inspection.

As was just pointed out, it can be quite useful to study the
SLE($\kappa, \rho$) 
in the context of conformal maps that fix $0$, $1$ and $\infty$.
It is therefore natural to define
$$
G_t (z) := \frac { g_t (z) - O_t }{W_t - O_t}
$$
where $(g_t, t \ge 0)$ is the Loewner chain associated with
SLE($\kappa, \rho$) and $O_t$ is the ``leftmost'' image of $0$ under
$g_t$.
The evolution equation for $G_t (z)$ is 
\begin {eqnarray*}
d_t G_t (z) 
&=& \frac {-\sqrt {\kappa}}{W_t - O_t}\, G_t (z) \,dB_t
\\
&&+
\frac {1}{(W_t - O_t)^2}\,  G_t (z)\, \Bigl(\frac {2}{G_t (z)
- 1} + \kappa - \rho - 2 \Bigr) \,dt
.\end {eqnarray*}   
If one then defines a time-change
$$u(t) = \int_1^t \frac {dv}{(W_v- O_v)^2}$$
and
$\tilde G_u (z) := G_{t(u)} (z)$ for all real $u$, 
then
$$
d\tilde G_u (z) 
= \tilde G_u (z) \sqrt {\kappa}\, d\tilde B_u 
+ \tilde G_u (z) \Bigl(\frac {2}{\tilde G_u(z) - 1}
+ \kappa - \rho - 2  \Bigr)\,dt,
$$
for a two-sided Brownian motion $\tilde B$
satisfying $d\tilde B_{u(t)}= -dB_t/(W_t-O_t)$.

\section{Equivalence of the frontiers of SLE$_6$ and Brownian motion}
\label{s.fileq}

 Brownian motion and
$SLE_6$ are both conformally invariant and
local.  We shall now see that this implies an fundamental equivalence between
the hulls that they generate.
 Some of the results presented in this section were announced in 
\cite {Wbarc}  
and have been presented in seminars for some years now.

\subsection {Full plane SLE$_6$ and planar Brownian motion}

The simplest version of the equivalence between the 
boundary of SLE$_6$ and planar Brownian motion
 involves full-plane $SLE_6$,
whose definition we now recall.
Let $\xi:\R\to\p\U$ be continuous.
It is well-known \cite{Pomm} that there is a unique one-parameter family
of conformal maps $f_t:\U\to\hat\C=\C\cup\{\infty\}$ such that the
inverses $g_t=f_t^{-1}$ satisfy Loewner's equation
\begin{equation}\label{e.radsle}
\p_t g_t(z) = - g_t(z)\,\frac{g_t(z)+\xi(t)}{g_t(z)-\xi(t)}\,,
\end{equation}
and the normalization (for all $t \in \R$)
\begin{equation}\label{e.radnormal}
\lim_{z\to\infty} z\,g_t(z)=e^t\,.
\end{equation}
Set $K_t:=\hat\C\setminus f_t(\U)$.
This is called the {\em hull} of the family
$(g_t)$.  Then $K_t\supset K_s$ when $t\ge s$, i.e.,
 $K_t$ is an increasing family of compact sets with
$\bigcap_t K_t = \{ 0 \}$. The relation (\ref{e.radnormal})
implies $\Cp(K_t) = t$, where $\Cp$ denotes
(logarithmic) capacity.  
(The capacity of a nonempty closed connected set $K\subset \C$
can be defined as $\lim_{z\to 0} \log |z\,g(z)|$,
where $g:\U\to\hat\C\setminus \fillplane(K)$ is any conformal map
satisfying $g(0)=\infty$.)
 
We remark that in the presence of~\eref{e.radsle}, if~\eref{e.radnormal}
holds for one $t\in\R$ it also holds for all $t\in\R$.
The proof of uniqueness is based on the fact that if $g_t$ and
$\hat g_t$ both satisfy \eref{e.radsle} and~\eref{e.radnormal}
for $t\in[-t_0,\infty)$ then $|g_t-\hat g_t|$ is necessarily
small away from $0$ if $t_0$ is large.  In other words, the far
away past matters very little.

Now let $\beta:\R\to\R$ be two-sided real Brownian motion with $\beta(0)=0$,
and let $b_0$ be random-uniform in $\p\U$ and independent
from $\beta$.  Set $\xi(t)= b_0\,\exp(i\, \sqrt \kappa  \beta(t))$.
With this choice of $\xi$, the above $(g_t,t\in\R)$ (or $(K_t, t \in \R)$)
is called
{\em full-plane SLE$_\kappa$}. As for other $SLE_\kappa$
(see \cite{RS,LSWlesl}), there
is a continuous path $\gamma:[-\infty,\infty]\to\hat\C$ with
$\gamma(-\infty) = 0$ such that $K_t=\fillplane(\gamma[-\infty,t])$.
Also,  given $\gamma[-\infty,s]$,  the
evolution of $\gamma(t), t\geq s$, is   the
same as the conformal image of the radial $SLE_\kappa$ path.  
In particular, full-plane $SLE_6$ satisfies the locality
property.

For the remainder of this section, we will fix $\kappa = 6$
and use $K_t$ to refer to the hull of $SLE_6$.
If $X_t, 0 \leq t < \infty$, 
is a complex Brownian motion starting at the origin, let
$\hat X_t=\fillplane(X[0,t])$ denote the {\em Brownian hull} at time $t$.
The {\em frontier} or {\em outer boundary} at time
$t$ is $\p \hat X_t$ and the elements
of $\bigcup_{t \geq 0} \p \hat X_t$
are called the {\em pioneer points}. 
 If $D$ is a simply
connected domain containing the origin, let
$\ttau_D = \inf\{t: K_t \in D\}$,
$\tau_D = \inf\{t: X_t \in D\}$. 

\begin{theorem}\label{fullplaneEQ}
Let $K_t$ denote the hull of full-plane $SLE_6$ and
$\hat X_t$ a planar Brownian hull as above.
Let $D\subset\C$ be a simply connected domain containing
the origin other than $\C$.
  Then $K_{\ttau_D}$
and $\hat X_{\tau_D}$ have the same law. 
\end{theorem} 
 
The proof of the theorem is based on the following lemma.

\begin{lemma}[Hitting measure for full-plane SLE$_6$]\label{slehit}
$K_{\ttau_D}\cap\p D$ is a single
point a.s.\ and the law of this point is the
same as the law of $X_{\tau_D}$, i.e., harmonic measure
from $0$.
\end{lemma}

\proofof{Theorem~\ref{fullplaneEQ} (using Lemma \ref{slehit})}
Consider some closed set $A\subset\hat \C$
such that $D'=D\setminus A$ is simply connected.
Then the law of $K_{\ttau_{D'}}\cap \p D'$ is harmonic measure from
$0$ on $\p D'$.  Consequently, the probability that
$K_{\ttau_{D'}}\cap A\ne \emptyset$ is the same as the 
probability that $X_{\tau_{D'}} \cap
\p D' \neq \emptyset.$    The former is equal to
$\P[K_{\ttau_D}\cap A\ne\emptyset]$ and the latter
is equal to $\P[\hat X_{\tau_D} \cap A\ne\emptyset]$.
Hence, the theorem follows by the corresponding analogue of
Lemma~\ref{l.uniq}.
\QED

\proofof{Lemma~\ref{slehit}}
Since radial $SLE_6$
satisfies locality~\cite{LSW2}, it follows that $\gamma(t)$ does too
and $K_{\ttau_D}$ is covariant with respect to conformal maps
$\phi:D\to D^*$.  (One uses here the fact mentioned above of stability
with respect to the far-away past.) 
Therefore, the law of $K_{\ttau_D}\cap\p D$ is harmonic measure.
\QED

\noindent {\bf Remark.}  Let $K^s = K_{T_{s \Disk}},$ $
\hat X^s = X_{\tau_{s \Disk}}$.  Then it is
easy to see that $ s^{-1} K^s,$ $s^{-1} \hat X^s$ are
continuous time Markov chains on the space of closed connected sets
$K$ contained
in $\overline \Disk$ with $0 \in K$  such that $K \cap \p \Disk$
is a single point.  This theorem can be interpreted
as saying that the two chains have the same invariant
distribution.  (By choosing full-plane $SLE_6$, we have
effectively started the chains in equilibrium.)
 However, it is not difficult (using
the fact that the Brownian motion crosses itself but
the $SLE_6$ path does not) to show
that the two chains are not the same.

\subsection {Chordal SLE$_6$ and reflected Brownian motion}

Let $A\subset \overline\H\setminus\{0\}$ be a closed subset of $\C$
such that the component of $0$ in $\overline\H\setminus A$ is
bounded.
Consider chordal $SLE_6$ in the upper-half plane and let 
$T:=\inf\{t: K_t\cap A\ne\emptyset\}$
denote the first time at which its hull $K_t\subset\overline\H$ hits $A$.
 
On the other hand, define a reflected Brownian motion in 
$\H$, $(B_t, t \ge 0)$ that is started from $0$ and reflected
on $\H$ with an angle 
$2\pi /3$
 pointing ``away'' from
the origin.
In other words,
$$
B_t = W_t + \int_0^t \bigl(e^{\pi i/3} d\ell^+_s + e^{2 \pi i/3} d\ell_s^-\bigr)
\,,
$$
where $W$ is standard
two-dimensional Brownian motion and $\ell^+$ (resp.\ $\ell^-$)
is a continuous process
that increases only when $B \in (0, \infty)$ (resp.\ $(-\infty, 0)$).
It is well-known (see, e.g., \cite {VW}) that this process
exists.
Let $\hat B_t:=\fillH(B[0,t])$, and
let $\tau=\tau_A$ denote the hitting time of $A$ by $B$.

\begin {theorem}\label{hpeq}
$K_T$ and $\hat B_\tau$ have the same law.
\end {theorem}

\proof
As in the proof of Theorem~\ref{fullplaneEQ},
it suffices to show that the law of
$K_T\cap A$ is the same as that of $\hat B_\tau\cap A$.

Let $S$ be the triangle with corners $0$, $e^{-2\pi i/3}$,
$e^{-\pi i/3}$.  Let $E_0$ denote the lower edge of $S$,
and let $E_1,E_2$ denote the other two edges.
Let $\phi:\H\setminus A\to S$ be
the conformal map which takes $0$ to $0$ and maps
$(\H \setminus A) \cap   A $ onto $E_0$.
The Cardy-Carleson formula for SLE$_6$ \cite{LSW1}
says that $\phi(K_T)\cap\p S$ is a uniformly chosen point
on $E_0$.  We need to prove the
same for $\phi(\hat B_\tau)\cap \p S$.
Reflected Brownian motion in smooth domains is conformally invariant up to 
a time change: the conformal invariance follows from
uniqueness and the fact that the conformal map
preserves angles up to the boundary by an application of
It\^o's formula.  (It\^o's formula is valid
at the reflection times too.)
Hence, it suffices to consider the hitting point on
$E_0$ by Brownian motion in $S$ starting from $0$,
reflected at angle 
$2 \pi /3$
 away from $0$ along the edges $E_1$ and $E_2$.

Let $n$ be large, and consider the triangular grid in $\overline S$
of mesh $1/n$, where the edges of $S$ are covered by edges of the grid.
Consider the Markov chain $Y$ on the vertices in $\overline S$
starting from $0$ with the following transition probabilities.
When $Y$ is at vertices interior to $S$ let $Y$ move with equal probability to
each of the neighbors.  At vertices 
$v\in (E_1\cup E_2)\setminus (E_0\cup\{0\})$
let $Y$ stay in $v$ with probability $1/6$, move to the neighbor below
$v$ on $E_1\cup E_2$ with probability $1/3$, and move to any of the
other three neighbors in $\overline S$ with probability $1/6$ each.
When $Y$ is at $0$ let it stay in $0$ with probability $1/3$
and move to each of its two neighbors in $\overline S$ with
probability $1/3$ each.  Let $Y$ stop when it hits $E_0$.
Induction shows that for $t=0,1,2,\dots$ conditioned on
$\Im\, Y(t)=h$ the distribution of $Y(t)$ is uniform among
vertices in $\overline S$ satisfying $\Im\, v =h$.
In particular, the vertex where $Y$ hits $E_0$ is uniform among
the vertices in $E_0$.  

It is not hard to verify that as $n$ tends
to $\infty$ the walk $Y(n^2 t)$ converges to the above
reflected Brownian motion in $S$.
\QED

Variations on this uniform 
hitting distribution property for reflected Brownian motions will be 
developed in \cite {DP}.
\medbreak

A {\em closed monotone class}\/ $\mathcal P$
is a collection of nonempty closed subsets $A\subset \hat \C$
that is closed in the Hausdorff topology and such that
$A\in\mathcal P$ and $A'\supset A$ implies $A'\in\mathcal P$
when $A'\subset\hat \C$ is closed.  For example,
$\mathcal P$ might be the collection of closed sets intersecting
some fixed closed set $Y$  or the collection of closed  connected
sets
whose capacity  is at least  $r$.
We now present a generalization of Theorems~\ref{fullplaneEQ}
and~\ref{hpeq}.

\begin{theorem}\label{monotonestop}
Let $\mathcal P$ be a closed monotone class.
Let $\ttau_{\mathcal P}:=\inf\{t\in\R: K_t\in\mathcal P\}$,
where $K_t$ is the hull of full-plane SLE$_6$
starting from $0$.
Let $X_t$ be planar Brownian motion starting from $X_0=0$ and
$\hat X_t:=\fillplane(X_t)$.
Let
$\tau_{\mathcal P}:=\inf\{t\ge 0: \hat X_t\in\mathcal P\}$.
Then $K_{\ttau_\mathcal P}$ and $\hat X_{\tau_\mathcal P}$
have the same law.

A corresponding generalization also holds for Theorem~\ref{hpeq}.
\end{theorem}

In the above, we take $K_\infty=\hat X_\infty=\hat \C$.
(This is relevant if $\ttau_{\mathcal P}=\infty$ or
$\tau_{\mathcal P}=\infty$ with positive
probability.)

\proof
As the proof in the chordal case is the same,
we will only treat the full-plane setting.
If $R > 0$, let $\mathcal P_R$ be the union of
$\mathcal P$ together with all closed sets intersecting
$R\,\p\U$.  Let $D$ be a simply connected domain containing $0$
other than $\C$.
Let $T_D =\sup\{t\in\R: K_t\subset D\}$ and
$\tau_D =\sup\{t\ge 0:\hat X_t\subset D\}$.
By Theorem~\ref{fullplaneEQ}, 
$$
\P[K_{T_D}\in\mathcal P_R]=\P[\hat X_{\tau_D}\in\mathcal P_R]\,.
$$
Observe that
$$
\P[K_{T_D}\in\mathcal P_R]=
\P[T_D\ge \ttau_{\mathcal P_R}]
=
\P[K_{\ttau_{\mathcal P_R}}\subset D]
+ \P[T_D=\ttau_{\mathcal P_R}]\,,
$$
and similarly
$$
\P[\hat X_{\tau_D}\in\mathcal P_R]=
\P[\hat X_{\tau_{\mathcal P_R}}\subset D] +
\P[\tau_{\mathcal P_R}=\tau_D]
\,.
$$
Shortly, we will prove
\begin{equation}\label{e.stptm}
\P[T_D=\ttau_{\mathcal P_R}]
=\P[\tau_{\mathcal P_R}=\tau_D]\,.
\end{equation}
Together with the above equalities this implies
$$
\P[K_{\ttau_{\mathcal P_R}}\subset D]=
\P[\hat X_{\tau_{\mathcal P_R}}\subset D] 
\,.
$$
The corresponding analogue of Lemma~\ref{l.uniq} then proves
that the laws of 
$K_{\ttau_{\mathcal P_R}}$ and $\hat X_{\tau_{\mathcal P_R}}$ are the same.
The theorem follows by letting $R\to\infty$.  It therefore remains
to prove~\eref{e.stptm}.

We claim that $\lim_{t\nearrow T_D} K_t = K_{T_D}$ a.s., in the Hausdorff
metric.  
By conformal invariance, it suffices to prove this for $D=\U$.
The times of discontinuity of $K_t$ (with respect to the
Hausdorff metric) are times $s$ where 
$K_s\setminus\lim_{t\nearrow s} K_t$ contains a nonempty open
set.  If $T_\U$ is a time of discontinuity with positive probability,
then the same would be true for $T_{r\U}$ for every $r>0$.
This would then contradict the fact that the
(expected) area of $K_{T_\U}$ is finite.
Hence $\lim_{t\nearrow T_D} K_t = K_{T_D}$ almost surely.
By monotonicity of $\mathcal P_R$, it follows that 
$\P[T_D=\ttau_{\mathcal P_R}]$ is equal to the probability
that $K_{T_D}\in\mathcal P_R$ but every compact subset of $K_{T_D}\cap D$
is not in $\mathcal P_R$.  Now~\eref{e.stptm} follows, because the
analogous argument applies to
$\P[\tau_{\mathcal P_R}=\tau_D]$
and the law of $\hat X_{\tau_D}$ is the same as that of
$K_{T_D}$.
\QED

\subsection {Chordal SLE$_6$ as Brownian motion reflected on its past hull}

By iterating the above results we will obtain an
\lq\lq emulation\rq\rq\ of chordal SLE$_6$ using reflected Brownian
motion.  Roughly, what we show is that the SLE$_6$ path is
Brownian motion that is reflected off its past filling with angle 
$2 \pi / 3$
towards infinity.

Let $(B_t^n, t \ge 0)_{n \ge 1}$
be a sequence of independent samples of 
reflected planar Brownian motion in $\H$ started from
$0$ that are reflected off the real axis with 
angle 
$ 2 \pi /3$
 away from $0$ (as before).
Define $\tilde B^1 = B^1$,  $\hat K^{0} := \emptyset$ and define inductively:
\begin {itemize}
\item $\tau_n:=\inf\{t\ge 0:|\tilde B^n_t-\tilde B^n_0|\ge \eps\}$,
\item $\hat K^n:=\fillH\bigl(\hat K^{n-1}\cup \tilde B^n[0,{\tau_n}]\bigr)$,
\item $\phi_n :\H\to\H \setminus \hat K^{n}$ is the conformal
 map normalized by $\phi_n(0)=\tilde B^{n}_{\tau_n}$, $\phi_n(\infty)= \infty$,
  $\phi_n'(\infty)=1$,
\item $\tilde B^{n+1}_t:=\phi_n\bigl(B^{n+1}_t\bigr)$.
\end {itemize}

\begin {corollary}
Fix $\eps>0$.
Let $\gamma$ denote the chordal SLE$_6$ path,
$K_t=\fillH\bigl(\gamma[0,t]\bigr)$ the SLE hull,
$T_0:=0$ and inductively,
$T_{n+1}:=\inf\{t\ge T_n:|\gamma(t)-\gamma(T_n)|\ge\eps\}$.
Then the sequence 
$(\hat K^0,\hat K^1, \ldots, )$ defined above has the same law
as the sequence $(K_{T_0}, K_{T_1}, \ldots )$.
Consequently, after reparameterization the path $\hat \gamma$
obtained by concatenating the paths
$\bigl(\tilde B^n_t, t\in[0,\tau_n]\bigr)$
stays within distance $2\eps$ from the path $\gamma$;
that is, $\sup_{t\ge 0} |\gamma(t)-\hat\gamma(t)|\le2\eps$.
\end{corollary}

\proof
This easily follows from Theorem~\ref{hpeq} and induction.
In the inductive step, we assume that
$(\hat K^{0},\dots,\hat K^n)$ and
$(K_{T_0},\dots,K_{T_n})$ have the same distribution.
Note that (when $n\ge 1$) this implies that
$(\hat K^{0},\dots,\hat K^n,\tilde B^n_{\tau_n})$  and
$(K_{T_0},\dots,K_{T_n}, \gamma(T_n))$
have the same distribution, since $\tilde B^n_{\tau_n}$
is a.s.\ the unique point in $\hat K^n$ at distance
$\eps$ from $\hat K^{n-1}$.
An application
of Theorem~\ref{hpeq} and conformal invariance now complete
the induction step and the proof.
\QED

One can also choose time  
sequences other than $(T_n, \tau_n)$. For instance, one can compare the 
sequence 
$(K_{t_n}, n \ge 1)$ (where $(t_n, n \ge 1)$ is a deterministic 
sequence) with 
$(\hat K^1, \hat K^2, \ldots  )$ 
where the definition of  
$\tau_n$ is replaced by
$$
\tau_n:=\inf\bigl\{t\ge 0: a( {\cal F}_\H ( \hat K^{n-1}
\cup \tilde B^n [0, t])) \ge 
2 t_n
\bigr\}.
$$
This time, one can for instance identify the tip of the SLE curve 
(or of the stopped reflected Brownian motion)
as the only accumulation point of cut points of the hull. We leave
the details to the interested reader.
 
\subsection {Non-equivalence of pioneer points and 
SLE$_6$}

Given these results, it is natural to try to better understand
the differences between planar Brownian motion 
and SLE$_6$.
How far does the equivalence go?
Consider, for example, the setting of Theorem~\ref{fullplaneEQ}.
Let $\theta(t):=\inf\{s\ge 0:\Cp(\hat X_s)\ge t\}$.
Theorem~\ref{monotonestop} implies that for
all $t\ge 0$, the distributions of
 $K_t$ and $\hat X_{\theta(t)}$
are the same. 
However, the processes $(K_t:t\ge 0)$
and $(\hat X_{\theta(t)}:t\ge 0)$ are different. 
In fact, it is not hard to show that
the joint distributions $(K_s,K_t)$ and
$(\hat X_{\theta(s)},\hat X_{\theta(t)})$ do not
agree when $s\ne t$. It is also true that the random
set  $X[0,\tau_\U]$
does not have the same distribution as
$\gamma[-\infty,T_\U)$; in fact, the first
has Hausdorff dimension $2$ and the latter
dimension $7/4$ \cite{Be2}.

The set $Z=\bigcup_{t\le\tau_\U} \p \hat X_t$
of   pioneer points of $X$ up to time $\tau_\U$
does have dimension $7/4$ a.s.\ \cite{LSWa}.
However, $Z$ does not have the same law as
$\gamma[-\infty,{\ttau_\U}] = \bigcup_{t \leq
\ttau_U} \p K_t$.  
We now give the outline of one 
possible proof of this fact.

\begin{figure}
\centerline{%
\SetLabels
\L(.52*.18)$z_0$\\
(.86*.35)$z_1$\\
(.65*.9)$z_3$\\
\R(.0*.5)$z_2$\\
\endSetLabels
\AffixLabels{\includegraphics*[height=2.3in]{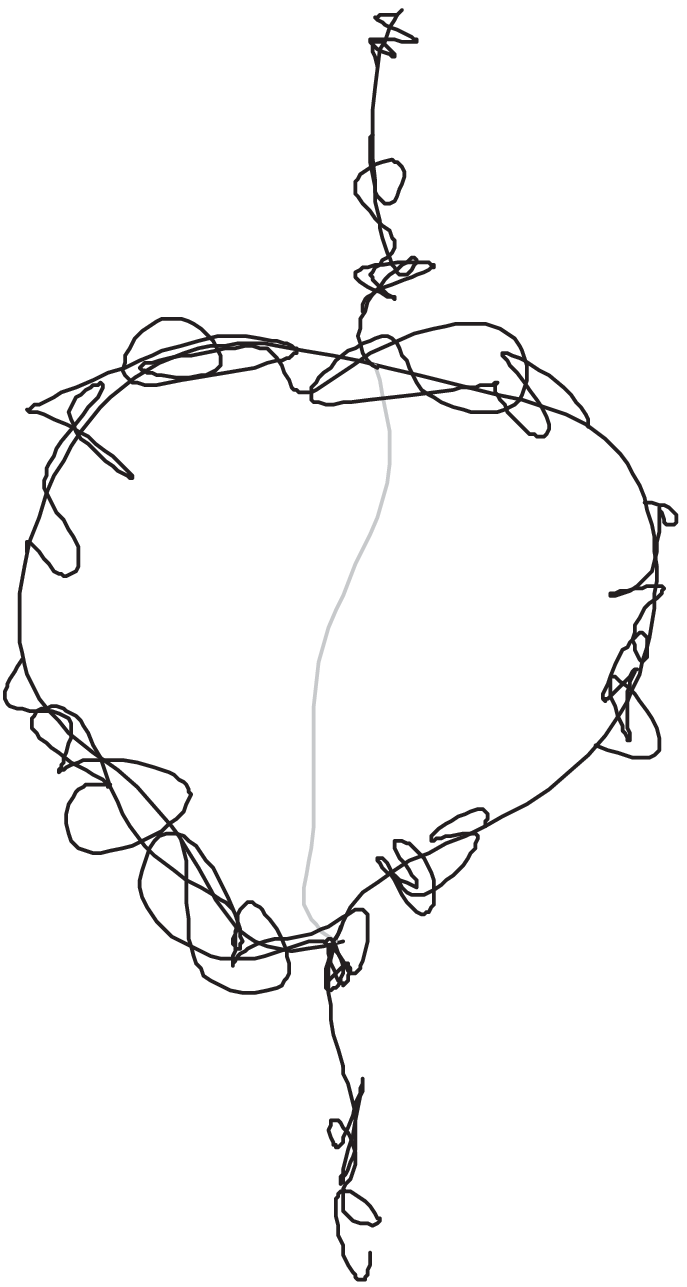}}%
}
\caption{\label{f.goodconfig}Good configurations for
the set of pioneer points.  The path indicated in gray represents
a segment of the Brownian path which does not contain any pioneer
points, since it was traversed after the loop surrounding it.}
\end{figure}

We say that $(z_0, z_1, z_2, z_3)$ is a good configuration
for $Z$ if:
\begin {itemize}
\item
$z_0$ and $z_3$ are cut-points of $Z$,
\item
any subpath of $Z$ from $z_0$ to $z_3$ goes through $z_1$ 
or $z_2$, and
\item 
there exist subpaths of $Z$ from $z_0$ to $z_3$ that 
go through $z_1$ (resp.\ through $z_2$) and not through $z_2$   
(resp.\ through $z_1$).
\end {itemize}
The set of good configurations of $Z$ is comparable to the 
set of cut-points of $Z$ in the sense that with 
positive probability, one can find four sets $Z_0, Z_1, 
Z_2, Z_3$ of Hausdorff dimension $3/4$ each (recall \cite {LSW2}
that the Hausdorff dimension of the set of cut points of the 
Brownian trace is $3/4$) such that any 
$(z_0, z_1, z_2, z_3) \in Z_0 \times Z_1\times Z_2\times Z_3$ 
is a good configuration.
This is due to the fact that for a  
Brownian path as shown in Figure~\ref{f.goodconfig}
$(z_0, z_1, z_2, z_3)$ is a good configuration.

On the other hand, a.s.\ the SLE$_6$ does not have good configurations.
For topological reasons ,  if $(z_0, z_1, z_2, z_3)$
is a good configuration for $\gamma[-\infty,\ttau_\U]$,
then at least one of the 
four points is a double point of $\gamma$
(hint: consider the two path-connected components of
$\gamma[-\infty,\ttau_\U]\setminus\{z_1,z_2\}$ which
contain $z_0$ and $z_3$, respectively and the
order in which $\gamma$ visits them).
In particular this point is 
simultaneously a local cut-point and a double point.
Such points do not exist for SLE$_6$ (as explained in
Remark~\ref{r.nocut}).
This argument can be made into a proof that the
the SLE$_6$ path image is not the same as the set of pioneer points of   
planar Brownian motion.
 
\subsection {Conditioned SLE$_6$}

We have seen
that the outer boundary of
a planar Brownian path looks (locally) like an SLE$_{8/3}$ path.
More precisely,
the right boundary (i.e., the right-hand side 
of the boundary of the filling) of a Brownian excursion  
in $\H$ is the path of SLE($8/3, 2/3$)
and the right boundary
of a reflected Brownian excursion with angle $3\pi/8$ 
is SLE($8/3$).
This gives some motivation to 
show that the outer boundaries of these Brownian excursions
have the same law as that of some conditioned SLE$_6$
processes, since this 
provides a description of  
the right boundary of conditioned SLE$_6$ 
in terms of variants of the SLE$_{8/3}$ paths.

In the spirit of the paper \cite {LSWup}, it is not 
difficult to prove that if one considers 
reflected Brownian motion $X$ in $\H$ (with any given reflection
angle), that is conditioned to hit (let $T$ be this
hitting time) 
$(-\infty, -1/\eps) \cup (1/\eps, \infty)$ before 
$(-1/\eps,  -\eps) \cup (\eps, 1/ \eps)$, and lets
$\eps \to 0$, the limiting law of $X[0,T]$
 is exactly that of a 
Brownian excursion (that does not touch the real line except
at the origin).
In particular, this implies that the filling of
chordal SLE$_6$ conditioned not to intersect the 
 real line (i.e., the limit when $\eps \to 0$ of $SLE_6$ 
conditioned not to intersect $(-1/\eps, - \eps) \cup (\eps, 1/\eps)$)
has law
$\P_1$. In particular, its right boundary is SLE($8/3,
 2/3$).

Similarly, the limit of the law of $X[0,T]$   
conditioned on $X[0,T] \cap (\eps, 1/\eps) = \emptyset$ is
simply the law of the reflected Brownian excursion.
If the reflection angle is 
$2\pi/3$ towards
infinity, as before, then the law of the right boundary of this
process is the 
SLE($8/3, -2/3$) path.
Hence, the right boundary of an SLE$_6$ conditioned not to 
intersect the positive half-line is exactly SLE($8/3, -2/3$).

\section {Remarks}
\label {rem}

Let us briefly sum up some of
the results that we have collected 
in the present paper 
concerning the description of the 
Brownian frontier.

\begin {itemize}
\item
The filling of the union of five independent excursions has the same law
as the filling of the union of eight independent chordal SLE$_{8/3}$.

\item
The right-boundary of a Brownian excursion from 0 to infinity in the 
upper half-plane reflected on the negative half-line with 
reflection angle $3 \pi / 8$ is SLE$_{8/3}$.
This law is symmetric with respect to reflection
 in the imaginary axis.  

\item
The right-boundary of a Brownian excursion is SLE($8/3, 2/3$).

\item
The right-boundary of the union of two independent Brownian
excursions is SLE($8/3, 2$).
As we shall mention shortly,
the right and left boundaries can also be viewed as non-intersecting
SLE$_{8/3}$'s.
 
\item
The Brownian bubble and the SLE$_{8/3}$ bubbles are identical
(up to scaling).

\end {itemize}

We conclude this paper by mentioning some closely related 
results that will be included in forthcoming papers:

\begin {itemize}
\item
Analogous problems in the
``radial case'', i.e.,   random subsets of the 
unit disk that contain one given boundary point and 
one given interior point, will be studied
in \cite {LSWrad}.  A radial restriction property
holds for SLE$_{8/3}$.  In particular, if $\gamma$
is a radial SLE$_{8/3}$ path in $\Disk$ from $1$ to
$0$, $A$ a compact set not containing $1$, 
such that $\Disk \setminus
A$ is simply connected and contains $0$, and
$\Psi_A$ is a conformal map from $\Disk \setminus
A$ onto $\Disk$ with $\Psi_A(0) = 0,$
  then
\[     \Prob\bigl[\gamma[0,\infty) \cap A = \emptyset\bigr]
    =  |\Psi_A'(0)|^{5/48} \; |\Psi_A'(1)|^{5/8} . \]

\item
In the spirit of   \cite {LW2},
  the Brownian half-space
intersection exponents
  computed in  \cite {LSW1,LSW3}
can be interpreted in terms of  non-intersection 
of independent sets defined under different
restriction measures. In particular, the measure $\P_2$
can be viewed as the filling of two SLE$_{8/3}$ 
that are conditioned not to intersect (of course, this 
event has probability $0$, so this has to be taken
as an appropriate limit).
 See e.g. \cite {W2}.

\item
In \cite {LSWloops},
a random countable set of loops $L$ in the plane called
the {\em Brownian loop soup}  is constructed.
Each $\gamma\in L$ is a loop, that is,
an equivalence class of
periodic, continuous maps from $\R$ to $\C$,
where $\gamma^1,\gamma^2$ are equivalent
if for some $r$, $\gamma^1(t) = \gamma^2(t+r)$
holds for all $t\in\R$.
Loosely speaking each $\gamma\in L$ is a {\em Brownian loop}. 
This loop soup is conformally invariant:
for any conformal map
$\Phi:D\to D'$ the sets
$\{ \Phi \circ \gamma: \gamma\in L,\, \gamma \subset D\}$ and
$\{ \gamma: \gamma\in L,\, \gamma \subset D'\}$ have
the same law, up to reparametrization of the loops.

It turns out that if one considers the set of loops in $\H$,
$L(\H)=\{\gamma\in L:\gamma\subset\H\}$,
and any Loewner chain $(K_t, t \in [0,T])$ generated
by a continuous curve, then 
another (equivalent) way to add a Poisson cloud of Brownian bubbles  
to the Loewner chain (as in our construction of the 
general restriction
$\P_\alpha$ measures) is  to add to the set $K_T$ all the loops
of $L(\H)$ that it intersects.
Therefore,
\begin {eqnarray*}
\lefteqn {
\P \bigl[ (K_T \cup \{ \hbox {cloud of bubbles} \} 
) \cap A = \emptyset
\bigr] }\\
&=&
\P\bigl [ \hbox {No loop in }L(\H)\hbox { intersects both } A
\hbox { and } K_T \bigr]
\\
&=&
\P\bigl [ K_T \cap (A \cup \{\hbox {loops that intersect A} \} )
= \emptyset \bigr].
\end {eqnarray*}
See \cite {LSWloops} for more details. This Brownian loop soup is then 
used in \cite {D2,W2}. 

\item
Restriction formulas can also be derived for SLE($\kappa,\rho$) processes,
and, combined with the loop soup, they shed 
light on the relation between SLE$_kappa$ and the 
outer boundary of SLE$_{16/\kappa}$ for $\kappa <4$. See \cite {D2}.

\end {itemize}

{\bigskip\noindent\bf Acknowledgments.}
We thank B\'alint Vir\'ag for fruitful conversations
and for permitting us to include Proposition~\ref{p.P1}.

\filbreak
\parindent=0pt

\vtop{  
\hsize=2.3in  
Gregory Lawler\\  
Department of Mathematics\\  
310 Malott Hall\\  
Cornell University\\  
Ithaca, NY 14853-4201, USA\\  
{lawler@math.cornell.edu}  
}  
\bigskip  
\vtop{  
\hsize=2.3in  
Oded Schramm\\  
Microsoft Corporation\\  
One Microsoft Way\\  
Redmond, WA 98052, USA\\  
{schramm@microsoft.com}  
}  
\bigskip  
\vtop{  
\hsize=2.3in  
Wendelin Werner\\  
D\'epartement de Math\'ematiques\\  
B\^at. 425\\  
Universit\'e Paris-Sud\\  
91405 ORSAY cedex, France\\  
{wendelin.werner@math.u-psud.fr}  
}  

\filbreak

\end {document}